\documentclass[aap,MSNbibl,seceqn,dvips]{arximspdf}
\usepackage{mathbh}
\usepackage{graphicx}

\doi{10.1214/12-AAP900} 
\volume{23}
\issue{6}
\pubyear{2013}
\firstpage{2259}
\lastpage{2289}

\makeatletter

\newcommand{\eqref}[1]{(\ref{#1})}
\newcommand{\ind}[1]{\mathbh{1}_{#1}}
\newcommand{\dint}{\mathrm{d}}

\newcommand{\eps}{\varepsilon}

\newtheorem{thmm}{Theorem}[section]
\newtheorem{corollary}[thmm]{Corollary}
\newtheorem{lemma}[thmm]{Lemma}
\newtheorem{proposition}[thmm]{Proposition}
\newtheorem{propositionn}{Proposition}[section]
\newproclaim{rem}[thmm]{Remark}
\newproclaim{assu}[thmm]{Assumption}
\makeatother

\begin{document}
\begin{frontmatter}

\title{Hitting time for Bessel processes---walk on moving spheres
algorithm (WoMS)\thanksref{T1}}
\runtitle{Hitting times for Bessel processes}
\thankstext{T1}{Supported in part by the French Research National Agency
(ANR) under the reference ANR-09-BLAN-0008-01.}

\begin{aug}
\author[a]{\fnms{Madalina} \snm{Deaconu}\ead[label=e1]{Madalina.Deaconu@inria.fr}}
\and
\author[b]{\fnms{Samuel} \snm{Herrmann}\corref{}\ead[label=e2]{Samuel.Herrmann@u-bourgogne.fr}}
\runauthor{M. Deaconu and S. Herrmann}
\affiliation{Inria and Universit{\'e} de Lorraine, and Universit{\'e}
de Bourgogne}


\address[a]{Institut Elie Cartan \\
\quad de Nancy (IECN)---UMR 7502\\
\quad and Project team TOSCA\\
INRIA Nancy Grand-Est\\
Universit{\'e} de Lorraine\\
B.P. 70239\\
54506 Vandoeuvre-l{\`e}s-Nancy Cedex\\
France\\
\printead{e1}}

\address[b]{Institut de Math{\'e}matiques\\
\quad de Bourgogne (IMB)---UMR 5584\\
Universit{\'e} de Bourgogne\\
B.P. 47 870\\
21078 Dijon Cedex\\
France\\
\printead{e2}}
\end{aug}

\received{\smonth{11} \syear{2011}}
\revised{\smonth{10} \syear{2012}}


\begin{abstract}
In this article we investigate the hitting time of some given
boundaries for Bessel processes. The main motivation comes from
mathematical finance when dealing with volatility models, but the
results can also be used in optimal control problems. The aim here is
to construct a new and efficient algorithm in order to approach this
hitting time. As an application we will consider the hitting time of a
given level for the Cox--Ingersoll--Ross process.
The main tools we use are on one side, an adaptation of the method of
images to this
particular situation and on the other side, the connection that exists
between Cox--Ingersoll--Ross processes and Bessel processes.
\end{abstract}

%


\end{frontmatter}

\section{Introduction}
The aim of this paper is to study the hitting time of some curved
boundaries for the Bessel process. Our main motivations come from
mathematical finance, optimal control and neuroscience. In finance
Cox--Ingersoll--Ross processes are widely used to model interest rates.
As an application, in this article we will consider the simulation of
the first hitting time of a given level for the CIR by using its
relation with the Bessel process.
In neuroscience the firing time of a neuron is usually modelled as the
hitting time of a stochastic process associated with the membrane
potential behavior; for introduction of noise in neuron systems, see
Part I Chapter 5 in \cite{gerstner}. The literature proposes different
continuous stochastic models like, for instance, the family of
integrate-and-fire models; see Chapter 10 in \cite{ermentrout}. Most of
them are related to the Ornstein--Uhlenbeck process which appears in a
natural way as extension of Stein's model, a classical discrete model.
In Feller's model, generalized Bessel processes appear as a more
realistic alternative to the Ornstein--Uhlenbeck process; see, for
instance, \cite{feller-mod} for a comparison of these models. Therefore
the interspike interval, which is interpreted as the first passage time
of the membrane potential through a given threshold is closely related
to the first hitting time of a curved boundary for some Bessel processes.

Our main results and the main algorithm are obtained for the case of
Bessel processes. We use in our numerical method the particular formula
that we obtain for the hitting time of some curved boundaries for the
Bessel process and the connection that exists between a Bessel process
and the Euclidean norm of a Brownian motion when calculating the
hitting position.
As an application we consider the hitting time of a given level for the
Cox--Ingersoll--Ross process. In order to obtain this, we use first of
all the connections that exist between CIR processes and Bessel
processes and second, the method of images for this particular situation.

The study of Bessel processes and their hitting times occupies a huge
portion of mathematical literature. Let us only mention few of them:
G{\"o}ing-Jaeschke and Yor~\cite{jaeschkeyor2003} consider a
particular case of CIR processes which are connected with radial
Ornstein--Uhlenbeck processes and their hitting times; L. Alili and P.
Patie \cite{alilipatie2010} investigate as a special situation the
Bessel processes via some boundary crossing identities for diffusions
having the time inversion property; recently, Salminen and  Yor
considered the hitting time of affine boundaries for the 3-dimensional
Bessel process \cite{salminenyor2011}.

In a recent paper Hamana and Matsumoto \cite{hamanamatsumoto2011}
gave explicit expressions for the distribution functions and the
densities of the first hitting time of a given level for the Bessel
process. Their results cover all the cases. Let us also mention a
recent work of Byczkowski, Malecki and Ryznar \cite
{byczkowskimaleckiryznar2011}. By using an integral formula for the
density of the first hitting time of a given level of the Bessel
process, they are able to obtain uniform estimates of the hitting time
density function.

In all these papers the formulae are explicit and are hard to use for a
numerical purposes as they exhibit Bessel functions. The main idea of
this paper is to get rid of this difficulty by using two important
tools: first of all the method of images that allow us to obtain, for
some particular boundaries, an explicit form for the density of the
hitting time, and second, the connection between $\delta$-dimensional
Bessel processes and the Euclidean norm of a $\delta$-dimensional
Brownian motion in order to get the simulated exit position. By
coupling these ingredients we are able to construct a numerical
algorithm that is easy to implement and very efficient and which
approaches the desired hitting time.

We will use here a modified version of the \textit{random walk on spheres}
method which was first introduced by Muller \cite{muller56} in 1956.
This procedure allows us to solve a Dirichlet boundary value problem.
The idea is to simulate iteratively, for the Brownian motion, the exit
position from the largest sphere included in the domain and centered in
the starting point. This exit position becomes the new starting point
and the procedure goes on until the exit point is close enough to the
boundary. Let us notice that the simulation of the exit time from a
sphere is numerically costly.

The method of images was introduced in 1969 by Daniels \cite
{daniels1969} as a tool to construct nonlinear boundaries for which
explicit calculations for the exit distribution for the Brownian motion
are possible. The method was developed also in Lerche~\cite
{lerche1986}. We adapt this method for the Bessel process by using the
explicit form of its density. For some particular curved boundaries we
can explicitly evaluate the density of the Bessel hitting time.

The paper is organized as follows. First we present some new results on
hitting times for Bessel processes. Second, we construct the new
algorithm for approaching the hitting time, the so called \textit{walk on
moving spheres algorithm}. Finally we present some numerical results
and as a particular application the evaluation of the hitting time of a
given level for the Cox--Ingersoll--Ross process.

\section{Hitting time for Bessel processes}

Bessel processes play an important role both in the study of
stochastic processes like Brownian motion and in various theoretical
and practical applications as, for example, in finance.

Let us consider the $\delta$-dimensional Bessel process starting from
$y$, the solution of the following stochastic differential equation:
%
\begin{equation}
\label{bessel-delta} \cases{ %
\displaystyle Z^{\delta,y}_t
= Z^{\delta,y}_0+\frac{\delta-1}{2}\int_0^t
\bigl(Z^{\delta,y}_s\bigr)^{-1}\, \dint s +
B_t,&\vspace*{2pt}
\cr
Z^{\delta,y}_0=y, \qquad y
\geq0,} %
\end{equation}
where $(B_t)_{t\geq 0}$ is a one-dimensional Brownian motion. We denote
%
\begin{equation}
\nu= \frac{\delta}{2}-1,
\end{equation}
the \emph{index} of this process. We call $\delta$ the \emph
{dimension} of the process. This terminology is coming from the fact
that, in the case of positive integer $\delta\in\mathbb{N}$, a $\delta
$-dimensional Bessel process can be represented as the Euclidean norm
of a $\delta$-dimensional Brownian motion. This will be a key point in
our numerical method.

The density of this process starting from $y$ is given by
%
\begin{equation}
\label{density}  p_y(t,x)=\frac{x}{t}\biggl(\frac{x}{y}
\biggr)^\nu\exp\biggl(- \frac{x^2+y^2}{2t}\biggr) I_\nu
\biggl(\frac{xy}{t}\biggr)\qquad \mbox{for } t>0, y>0, x\geq0,\hspace*{-35pt}
\end{equation}
where $I_\nu(z) $ is the Bessel function whose expression gives
%
\begin{equation}
\label{bessel-function} I_\nu(z) = \sum_{n=0}^\infty
\biggl( \frac{z}{2} \biggr) ^{\nu+2n} \frac{1}{n ! \Gamma(\nu+n+1)}.
\end{equation}
When starting from $y=0$, the density of $Z^{\delta,0}_t$ is
%
\begin{equation}
\label{density-bessel-0} p_0(t,x) = \frac{1}{2^{\nu}}
\frac{1}{t^{\nu+1}}\frac{1}{\Gamma
(\nu+1)} x^{\delta-1}\exp \biggl(-
\frac{x^2}{2t} \biggr)\qquad\mbox{for } t>0, x\geq0.
\end{equation}

\subsection{The method of images for Bessel processes}
In this section, we investigate the first hitting time of a curved
boundary for the Bessel process starting from the origin. Let $\psi(t)$
denote the boundary, and introduce the following stopping time:
%
\begin{equation}
\label{taupsi} \tau_\psi=\inf\bigl\{ t> 0; Z^{\delta,0}_t
\geq\psi(t)\bigr\}.
\end{equation}
For some suitable choice of the boundary, the distribution of $\tau_\psi
$ can be explicitly computed. The idea is given by the following remark
on the method of images (see, e.g., \cite{daniels1969} for the
origin of this method and \cite{lerche1986} for a complete presentation):

\textit{Fundamental idea}.
Suppose that $F$ is a positive $\sigma$-finite measure satisfying some
integrability conditions (to be specified later on), and define
%
\begin{equation}
u(t,x)=p_0(t,x)-\frac{1}{a}\int_{\mathbb{R}_+}p_y(t,x)
F(\dint y)
\end{equation}
for some real constant $a>0$. Let
\[
\psi(t) =\inf\bigl\{ x\in\mathbb{R}; u(t,x) < 0\bigr\}\qquad \mbox{ for all } t>0.
\]
Then $u(t,x)$ is solution of the partial differential equation
%
\begin{equation}
\label{PDEbessel} \cases{ %
\displaystyle\frac{\partial u}{\partial t}(t,x)
= \frac{1}{2} \frac
{\partial^2 u}{\partial x^2}(t,x) -\frac{\delta-1}{2}
\frac
{\partial}{\partial x} \biggl(\frac{1}{x} u(t,x) \biggr), &\quad $\mbox{on }
\mathbb{R}_+\times\mathbb{R},$
\vspace*{2pt}\cr
u\bigl(t,\psi(t)\bigr) = 0,& \quad$\mbox{for all } t>0,$
\vspace*{2pt}\cr
u(0,\cdot)=\delta_0(\cdot), & \quad$\mbox{on } (-\infty, \psi(0+)].$}\hspace*{-35pt}
\end{equation}
From this remark we deduce an interesting expression for the hitting
time. We can prove that
\[
\tau_\psi=\inf\bigl\{t>0; u\bigl(t, Z^{\delta,0}_t
\bigr)=0\bigr\}.
\]
This means simply that in order to obtain information on the hitting
time it suffices to look for $u(t, Z^{\delta, 0}_t)=0$.

Let us express this in a general result.
%
\begin{thmm}
\label{thmgeneralsetting}
Let $F(\dint y)$ be a positive $\sigma$-finite measure such that $\int_0^\infty p_0(t,\break \sqrt{\varepsilon} y) F(\dint y) <\infty$ for all
$\varepsilon>0$. Let $a>0$ and define the function
%
\begin{equation}
\label{generalu} u(t,x)=p_0(t,x)-\frac{1}{a}\int
_{\mathbb{R}_+}p_y(t,x) F(\dint y).
\end{equation}
Consider $\psi(t)$ such that $u(t,\psi(t))=0$. Then the probability
density function of $\tau_\psi$ is given by
%
\begin{eqnarray}
\label{distributiontau} \qquad&&\mathbb{P}_0(\tau_\psi\in\dint t)
\nonumber
\\[-8pt]
\\[-8pt]
\nonumber
&&\qquad= \biggl[-\frac{1}{2} \frac
{\partial u}{\partial x}(t,x)\Big|
_{x=\psi(t)}+\frac{1}{2} \frac{\partial u}{\partial x}(t,x)\Big|_{x=0} - \frac
{\delta-1}{2x} u(t,x)\Big|_{x=0} \biggr]\,
\dint t.
\end{eqnarray}
\end{thmm}
%
%
\begin{pf} We will only point out the main ideas for the proof in
this case as it follows mainly the ideas introduced in \cite
{lerche1986}. A complete description of the method and this result for
the Brownian motion case can be found in~\cite{lerche1986}.

Let us consider
%
\begin{equation}
\label{function-u-1} u(t,x)= p_0(t,x)-\frac{1}{a} \int
_{\mathbb{R_+}} p_y(t,x)F(\dint y),
\end{equation}
where $F(\dint y)$ is a measure on $\mathbb{R}_+$.
We consider $\psi(t)$ the solution of $u(t,\psi(t))=0$.
Let us define as before $\tau_\psi= \inf\{ t\geq0; Z^{\delta,0}_t\geq
\psi(t)\}$. Then $u(t,x)\,\dint x=\mathbb{P}(\tau_\psi>t, Z^{\delta
,0}_t\in\dint x)$ and
%
\begin{equation}
\label{probahittingtime} \mathbb{P}_0(\tau_\psi> t) =\int
_0^{\psi(t)} u(t,x)\,\dint x.
\end{equation}
In order to get the distribution of $\tau_\psi$ we have to evaluate the
derivative of \mbox{$\mathbb{P}_0(\tau_\psi> t )$}. By using equality (\ref
{probahittingtime}) we obtain
%
\begin{eqnarray}
&& \mathbb{P}_0(\tau_\psi\in
\dint t)\nonumber\\
&&\qquad= \biggl(-\psi'(t) u\bigl(t, \psi(t)\bigr) -\int
_0^{\psi(t)} \frac{\partial u}{\partial t} (t,x) \,\dint x \biggr)\,
\dint t
\\
&&\qquad= \biggl(-\frac{1}{2}\int_0^{\psi(t)}
\frac{\partial^2
u}{\partial x^2}(t,x)\,\dint x +\frac{\delta-1}{2}\int_0^{\psi(t)}
\frac{\partial}{\partial x} \biggl(\frac{1}{x} u(t,x) \biggr)\,\dint x \biggr)\,\dint
t,\nonumber
\end{eqnarray}
as $u(t,x)$ is solution of partial differential equation (\ref{PDEbessel}).
We thus obtain
%
\begin{eqnarray}
\label{probatau} %
&& \mathbb{P}_0(
\tau_\psi\in\dint t)  = \biggl( -\frac{1}{2}
\frac{\partial u}{\partial x}(t,x)\Big| _{x=\psi(t)}+\frac{\delta
-1}{2\psi(t)} u\bigl(t,\psi(t)
\bigr)
\nonumber
\\[-8pt]
\\[-8pt]
\nonumber
&&\hspace*{83pt}{}+  \biggl(\frac{1}{2} \frac{\partial u}{\partial x}(t,x) -
\frac{\delta-1}{2x} u(t,x) \biggr)\Big|_{x=0} \biggr)\,\dint t.
\end{eqnarray}
As $\frac{\delta-1}{2\psi(t)} u(t,\psi(t))=0$, and this ends the proof
of the theorem.
\end{pf}
The idea behind the method of images is that for some particular forms
of $F(\dint y),$ we can derive explicit formulae of the hitting time
distribution. More precisely:
%
\begin{proposition}
\label{propfirstboundary}
Let us denote, for $\delta=2\nu+2 > 0$ and
$a > 0$ by
\[
\operatorname{Supp} (\tau_\psi)= \biggl[ 0, \biggl(\frac{a}{\Gamma(\nu+1)2^\nu}
\biggr)^{{1}/{(\nu+1)}} \biggr].
\]
We define, for all $t\in \operatorname{Supp} (\tau_\psi)$, the function
%
\begin{equation}
\label{eqajout} \psi_a(t)=\sqrt{2t\log\frac{a}{\Gamma(\nu+1) t^{\nu+1}2^\nu}}.
\end{equation}
Then the probability density of $\tau_\psi$ has its support in $\operatorname{Supp}
(\tau_\psi)$ and is given by
%
\begin{equation}
\label{densitetau} \mathbb{P}_0(\tau_\psi\in\dint t) =
\frac{1}{2at} \biggl( 2t \log \frac{a}{\Gamma(\nu+1)t^{\nu+1}2^\nu} \biggr)^ {\nu+1}
\ind{\operatorname{Supp}(\tau _\psi)}(t)\,\dint t.
\end{equation}
\end{proposition}
\begin{pf}
By using the expression in (\ref{density}) we remark first that
%
\begin{equation}
\label{change-x-y} y^{2\nu+1} p_y(t,x)= x^{2\nu+1}
p_x(t,y).
\end{equation}
Let us consider, as in Theorem \ref{thmgeneralsetting},
%
\begin{equation}
\label{function-u-2} u(t,x)= p_0(t,x)-\frac{1}{a} \int
_{\mathbb{R_+}} p_y(t,x)F(\dint y),
\end{equation}
with $F(\dint y) = y^{2\nu+1} \ind{\{y >0\}}\,\dint y$. In this situation
the function $u$ defined in \eqref{function-u-2} gives
%
\begin{eqnarray}
\label{specific-F} %
u(t,x) &=& p_0(t,x)
-\frac{1}{a} x^{2\nu+1}
\nonumber
\\[-8pt]
\\[-8pt]
\nonumber
& =& \biggl( \frac{1}{2^{\nu}}\frac{1}{t^{\nu+1}}\frac{1}{\Gamma
(\nu+1)}\exp \biggl(-
\frac{x^2}{2t} \biggr)-\frac{1}{a} \biggr)x^{2\nu+1}.
\end{eqnarray}
For simplicity we will write $\psi$ instead of $\psi_a$. Following the
result in the Theorem~\ref{thmgeneralsetting}, we are looking for
$\psi(t)$ such that $u(t,\psi(t))=0$. This yields
%
\begin{equation}
\label{curve} x=\psi(t)=\sqrt{2t\log\frac{a}{\Gamma(\nu+1) t^{\nu+1}2^\nu}}
\end{equation}
under the obvious condition $t^{\nu+1} \leq\frac{a}{\Gamma(\nu+1) 2^\nu
}$.

We can now notice that
\[
p_0\bigl(t,\psi(t)\bigr) = \frac{1}{a}\bigl(\psi(t)
\bigr)^{2\nu+1},
\]
and we can prove easily that
\[
\frac{\partial u}{\partial x} (t,x)=(\delta-1)\frac{u(t,x)}{x} -\frac{x}{t}p_0(t,x).
\]
We obtain, after replacing in (\ref{distributiontau}) and after
applying the Theorem \ref{thmgeneralsetting}, for this particular case,
\begin{eqnarray*}
 \mathbb{P}_0(\tau_\psi\in
\dint t)&=&\frac{1}{2t} \psi(t) p_0\bigl(t,\psi(t)\bigr)\,\dint t
\\
& = &\frac{1}{2at} \psi^{2\nu+2} (t)\,\dint t
\\
&= &\frac{1}{2at} \biggl( 2t \log\frac{a}{\Gamma(\nu+1)t^{\nu
+1}2^\nu} \biggr)^{\nu+1}\,
\dint t,
\end{eqnarray*}
which gives the desired result.
\end{pf}
The second boundary which allows us to express explicit results is
obtained by using the Markov property for the Bessel process.
%
\begin{proposition}
\label{propsecondboundary}
Let us, for $\delta=2\nu+2 > 0$, $s>0$ and $a >0$ fixed, denote by
\[
\operatorname{Supp}(\tau_\psi) = \cases{ %
[0, +\infty), &\quad
$\mbox{for }  a\geq1,$
\vspace*{2pt}\cr
\displaystyle\biggl[ 0, \frac{s}{( {1}/{a})^{{1}/{(\nu+1)}}-1} \biggr], & \quad$\mbox{for }
0<a<1.$}
\]
For $t\in \operatorname{Supp}(\tau_\psi)$ we define the function
%
\begin{equation}
\label{secondboundary} \psi_a(t)=\sqrt{\frac{2t(t+s)}{s} \biggl[ (
\nu+1)\log \biggl(1+\frac
{s}{t} \biggr) +\log a \biggr]}.
\end{equation}
%
Then the probability density function of the hitting time $\tau_\psi$
is given by
%
\begin{eqnarray}
\label{probatausecondboundary} %
&&\mathbb{P}_0(\tau_\psi\in\dint t)\nonumber\\
&&\qquad=
\frac{1}{\Gamma(\nu+1)}\frac{1}{t} \biggl(\frac
{t+s}{s}
\biggr)^{\nu} \biggl[\log \biggl( a \biggl(\frac{t+s}{t}
\biggr)^{\nu
+1} \biggr) \biggr]^{\nu+1} \\
&&\qquad\quad{}\times\exp \biggl[-
\frac{t+s}{s}\log \biggl(a \biggl(\frac{t+s}{t} \biggr)^{\nu+1}
\biggr) \biggr]\ind{\operatorname{Supp}(\tau_\psi)}(t) \,\dint
t.\nonumber
\end{eqnarray}
\end{proposition}
\begin{pf} We will only sketch the proof as it follows the same
ideas as
the one of the Theorem \ref{thmgeneralsetting}. Let us consider the measure
$F(\dint y)= p_0(s,y)\,\dint y$ for $s>0$ fixed.
Then, when evaluating the corresponding $u(t,x)$, we have
\begin{eqnarray*}
u(t,x) 
& =& p_0(t,x)-
\frac{1}{a}\int_{\mathbb{R}_+} p_0(s,y)
p_y(t,x)\,\dint y
\\
& =& p_0(t,x)-\frac{1}{a} p_0(t+s,x)
\\
& =& \frac{1}{2^\nu}\frac{1}{\Gamma(\nu+1)}x^{2\nu+1} \biggl[
\frac{1}{t^{\nu+1}} \exp \biggl( -\frac{x^2}{2t} \biggr) -\frac
{1}{a}
\frac{1}{(t+s)^{\nu+1}} \exp \biggl(-\frac
{x^2}{2(t+s)} \biggr) \biggr],
\end{eqnarray*}
by using the Markov property.
We obtain the form of $\psi(t)$ by the condition $u(t,\psi(t))=0$
which gives
%
\begin{eqnarray}
\psi(t) = \sqrt{\frac{2t(t+s)}{s} \biggl[ (\nu+1)\log \biggl(1+
\frac
{s}{t} \biggr) +\log a \biggr]},
\nonumber
\\[-8pt]
\\[-8pt]
 \eqntext{\cases{ %
\mbox{for } t\geq0 \mbox{ if } a\geq1,
\vspace*{2pt}\cr
\mbox{for } t\leq\displaystyle\frac{s}{({1}/{a})^{{1}/{(\nu+1)}}-1} \mbox{ if } a <
1.}}
\end{eqnarray}
In order to obtain the distribution of $\tau_\psi$, one has only to evaluate
%
\begin{equation}
\frac{\partial u}{\partial x} (t,x) = (\delta-1)\frac{u(t,x)}{x} -\frac{s}{t(t+s)} x
p_0(t,x),
\end{equation}
and $\frac{u(t,x)}{x}$ for $x=0$ and $x=\psi(t)$ and replace the
values in the general form (\ref{probatau}). The expression (\ref
{probatausecondboundary}) follows.
\end{pf}
%
\begin{rem}
We can notice that the function $\psi_a(t)$ defined by (\ref
{secondboundary}) satisfies, for large times,
\[
\cases{ %
\psi_a(t)\simeq \sqrt{t}, &\quad$\mbox{for } a =1,$
\vspace*{2pt}\cr
\psi_a(t)\simeq t, &\quad  $\mbox{for all }a >1.$}
\]
In particular, we can approach large times by considering this kind of boundary.
\end{rem}
A new boundary can be obtained by using the Laplace transform of the
square of the $\delta $-dimensional Bessel process starting from $x$.
More precisely:
%
\begin{proposition}
\label{propthirdboundary}
Let, for $\delta=2\nu+2 > 0$, $\lambda>0$ and $a>0$ fixed,
%
\begin{equation}
\label{psilaplace} \psi_a(t)=\frac{\lambda t}{1+2\lambda t}+ t\sqrt{ \biggl(
\frac{\lambda}{1+2\lambda t} \biggr)^2 +\frac{2}{t} \log
\frac{a(1+2\lambda t)^{\nu+1}}{2^\nu t^{\nu+1}\Gamma(\nu+1)}}
\end{equation}
for all $t\in \operatorname{Supp}(\tau_\psi)$, where $\operatorname{Supp}(\tau_\psi)$ is defined by
%
\begin{equation}
\operatorname{Supp}(\tau_\psi) = \cases{ %
\displaystyle \biggl[0,
\frac{1}{(2^\nu\Gamma(\nu+1)/a)^{1/(\nu+1)}-2\lambda} \biggr], \vspace*{2pt}\cr
\qquad \displaystyle\mbox{if } \lambda< \frac{1}{2} \biggl(
\frac{2^\nu\Gamma(\nu+1)}{a} \biggr)^{{1}/{(\nu+1)}},
\vspace*{2pt}\cr
[0, +\infty),\vspace*{2pt}\cr
 \qquad \displaystyle\mbox{if }
\lambda\geq\frac{1}{2} \biggl(\frac
{2^\nu\Gamma(\nu+1)}{a} \biggr)^{{1}/{(\nu+1)}}.}
\end{equation}
Then the probability density function of the hitting time is given by
%
\begin{eqnarray}
\label{probatauthirdboundary} &&\mathbb{P}_0(\tau_\psi\in
\dint t)
\nonumber
\\[-8pt]
\\[-8pt]
\nonumber
&&\qquad=\sqrt{ \biggl(\frac{\lambda
}{1+2\lambda t} \biggr)^2 +
\frac{2}{t} \log\frac{a(1+2\lambda
t)^{\nu+1}}{2^\nu t^{\nu+1}\Gamma(\nu+1)}}p_0\bigl(t,\psi(t)\bigr)\ind
{\operatorname{Supp}(\tau_\psi)}(t)\,\dint t.\hspace*{-35pt}
\end{eqnarray}
\end{proposition}
\begin{pf}
We present only the main ideas as the result follows as above from the
general method in Theorem \ref{thmgeneralsetting} applied to the
measure $F(\dint y) = y^{2\nu+1} e ^{-\lambda y^2}\ind{\{y\geq0\}
}\,\dint y$. For this measure $u(t,x)$ takes the form
\begin{eqnarray*}
 u(t,x) & = & p_0(t,x)-
\frac{1}{a}\int_{\mathbb{R}_+} p_y(t,x) F(\dint y)
\\
&=& p_0(t,x) -\frac{1}{a}\int_{\mathbb{R}_+}
p_y(t,x)y^{2\nu+1} e ^{-\lambda y^2}\,\dint y
\\
&=& p_0(t,x) -\frac{1}{a}\int_{\mathbb{R}_+}x^{2\nu+1}p_x(t,y)
e ^{-\lambda y^2}\,\dint y
\\
&=& p_0(t,x) -\frac{x^{2\nu+1}}{a} \mathbb{E} \bigl(e^{-\lambda
Z_t^{\delta,x}}
\bigr).
\end{eqnarray*}
By using the expression of the Laplace transform for $Z_t^{\delta,x}$
we obtain
%
\begin{equation}
\label{TLboundary} u(t,x) = p_0(t,x) -\frac{x^{2\nu+1} }{a}
\frac{1}{(1+2\lambda
t)^{\delta/2}}\exp \biggl(-\frac{\lambda x}{1+2\lambda t} \biggr).
\end{equation}
We consider first the equality $u(t,\psi(t))=0$ in (\ref
{TLboundary}), and this gives the form of $\psi(t)$ in (\ref{psilaplace}).
Afterwards, we can evaluate once again in this particular situation
\[
\frac{\partial u}{\partial x}(t,x) = (\delta-1)\frac{u(t,x)}{x} - \biggl(
\frac{x}{t} -\frac{\lambda t}{1+2\lambda t} \biggr) p_0(t,x).
\]
For this particular case, there is only one nonvanishing term in
expression (\ref{distributiontau}) of
$\mathbb{P}_0(\tau_\psi \in \dint t)$, that is, the term $-
(\frac{x}{t} -\frac{\lambda t}{1+2\lambda t} ) p_0(t,x)$ of
$\frac{\partial u}{\partial x}(t,x)$ for $x=\psi(t)$, and this is
exactly given by the right-hand side of formula (\ref
{probatauthirdboundary}).
\end{pf}
%
\begin{corollary}
The previous results give, for $\delta=2$:
\begin{longlist}[(1)]
\item[(1)] for $a>0$, $0\leq t\leq a$ and $\psi(t)=\sqrt{2t\log
\frac{a}{t}}$, the density of the hitting time $\tau_\psi$ is
\[
\label{densitetau1} \mathbb{P}_0(\tau_\psi\in\dint t) =
\frac{1}{2a} \log \biggl( \frac{a}{t} \biggr)\ind{\{0\leq t\leq a
\}}(t)\,\dint t;
\]
\item[(2)] for $s>0$, $a>0$, $0\leq t\leq\frac{sa}{1-a}$ and $\psi
(t) = \sqrt{\frac{2t(t+s)}{s}\log (a\frac{t+s}{t} )}$,
the probability density function of $\tau_\psi$ is given by
\begin{eqnarray*}
\label{densitetau2} &&\mathbb{P}_0(\tau_\psi\in\dint t)\\
&&\qquad =
\frac{t+s}{t} \log \biggl( a \frac{t+s}{t} \biggr) \exp \biggl[-
\frac{t+s}{t}\log \biggl( a\frac
{t+s}{t} \biggr) \biggr]\ind{ \biggl\{0
\leq t\leq{sa}/{(1-a)} \biggr\}}(t) \,\dint t;
\end{eqnarray*}
\item[(3)] for $a>0$ and $\psi(t) =\frac{\lambda t}{1+2\lambda
t}+t\sqrt{ (\frac{\lambda}{1+2\lambda t} )^2 +\frac
{2}{t}\log\frac{a(1+2\lambda t)}{t}}$, for $t\in\break  \operatorname{Supp}(\tau_\psi)$, where
%
\begin{equation}
\operatorname{Supp}(\tau_\psi) = \cases{ %
[0, +\infty), &\quad
$\displaystyle\mbox{if } \lambda\geq {\frac{1}{2a}},$
\vspace*{2pt}\cr
\displaystyle\biggl[ 0, \frac{a}{1-2\lambda a} \biggr], &\quad $\displaystyle\mbox{if } \lambda<
\frac{1}{2a},$}
\end{equation}
the probability density function of $\tau_\psi$ is
\[
\label{densitetau3} \mathbb{P}_0(\tau_\psi\in\dint t) =
\sqrt{ \biggl(\frac{\lambda
}{1+2\lambda t} \biggr)^2 +\frac{2}{t}\log
\frac{a(1+2\lambda
t)}{t}} p_0\bigl(t, \psi(t)\bigr)\ind{\operatorname{Supp}(
\tau_\psi)}(t)\,\dint t.
\]
\end{longlist}

\end{corollary}

\subsection{Approximation of the first hitting time for Bessel
processes starting from the origin}
In this section we will construct a stepwise procedure, the so-called
\emph{random walk on moving spheres (WoMS)} algorithm, which allows us
to approach the first time the standard Bessel process hits a given
level $l>0$. Of course, this stopping time $\tau_l=\inf\{t> 0;
Z^{\delta,x}_t=l\}$ can be characterized by its well-known Laplace
transform computed by solving an eigenvalue problem. Indeed if
$(Z^{\delta,x}_t, t\ge0)$ is the Bessel process starting from $x$, of
index $\nu=\frac{\delta}{2}-1$, then for $\nu>0$ and $x\le l$, we get
\[
\mathbb{E}_x\bigl[e^{-\lambda\tau_l}\bigr]=\frac{x^{-\nu}I_\nu(x\sqrt{2\lambda
})}{l^{-\nu}I_\nu(l\sqrt{2\lambda})}, x>0\quad
\mbox{and} \quad\mathbb {E}_0 \bigl[ e^{-\lambda\tau_l} \bigr]=
\frac{(l\sqrt{2\lambda})^\nu
}{2^\nu\Gamma(\nu+1)}\frac{1}{I_\nu(l\sqrt{2\lambda})}.
\]
Here $I_\nu$ denotes the modified Bessel function. This Laplace
transform can be used to describe the following tail distribution:
Ciesielski and Taylor \cite{Ciesielski-Taylor} proved that, for $\delta
=2\nu+2\in\mathbb{N}^*$,
\[
\mathbb{P}_0(\tau_l>t)=\frac{1}{2^{\nu-1}\Gamma(\nu+1)} \sum
_{k=1}^{\infty}\frac{j_{\nu,k}^{\nu-1}}{\mathcal{J}_{\nu+1}(j_{\nu
,k})} e^{(-{j_{\nu,k}^2}/{(2l^2)})t},
\]
where $\mathcal{J}_\nu$ is the Bessel function of the first kind, and
$(j_{\nu,k})_{\nu,k}$ is the associated sequence of its positive
zeros.

We are looking for a numerical approach for the hitting time and these
formulae are not easy to handle and approach, in particular we cannot
compute the inverse cumulative function! The aim of this section is to
construct an easy and efficient algorithm without need of inverting the
Laplace transform and without using discretization schemes since the
hitting times are unbounded. In the next step we will extend this
procedure to the hitting time of time-dependent boundaries like
straight lines, useful in the description of the hitting time of a
given level for the CIR process (the Laplace transform is then unknown).
%
%
\subsubsection{\texorpdfstring{Hitting time of a given level for the Bessel process
with positive integer dimension $\delta$}
{Hitting time of a given level for the Bessel process
with positive integer dimension delta}}
\label{sect0-horiz}
Let us consider $\delta$ independent one-dimensional Brownian motions
$(B^{(k)}_t,t\ge0)$, $1\le k\le\delta$. Then the Bessel process of
index $\delta$ starting from 0, satisfies the following property:
\[
\bigl(Z^{\delta,0}_t, t\ge0\bigr) \mbox{ has the same
distribution as } \bigl(\sqrt{ \bigl(B^{(1)}_t
\bigr)^2+\cdots+ \bigl(B^{(\delta)}_t
\bigr)^2}, t\ge0 \bigr).
\]
Let
%
\begin{equation}
\label{tau-l} \tau_l = \inf\bigl\{ t\geq0; Z^{\delta,0}_t
\geq l\bigr\}.
\end{equation}
In particular, we can express $\tau_l$ by using the first time when the
$\delta$-dimensional Brownian motion $\mathbf{B}=(B^{(1)},\ldots,B^{(\delta
)})$ exits from the Euclidean ball $D$ centered in the origin with
radius $l$. Approximating the exit time and the exit position for the
2-dimensional Brownian motion of a convex domain was already studied by
Milstein \cite{Milstein-97}. He used the so-called random walk on
spheres algorithm which allows one to approach the exit location and
the exit time through an efficient algorithm. The exit position is
really simple to obtain (as it is uniformly distributed on the circle)
while the exit time is much more difficult to approach. That is why we
will construct an adaptation of this initial procedure in order to
obtain nice and efficient results concerning the Bessel process exit
time.

Let us introduce now our \emph{walk on moving spheres} ($\mathit{WoMS}$)
algorithm. We first define a continuous function $\rho\dvtx \mathbb{R}^2\to
\mathbb{R}_+$ which represents the distance to the boundary of $D$,
%
\begin{equation}
\label{defderho} \rho(x)=\inf\bigl\{\Vert x-y\Vert; y\in D^{c}\bigr
\}=l-\Vert x\Vert.
\end{equation}
For any small enough parameter $\varepsilon>0$, we will denote by
$D^\varepsilon$ the sphere centered at the origin with radius
$l-\varepsilon$,
%
\begin{equation}
\label{defDepsilon} D^{\varepsilon} = \{ x\in D; \Vert x\Vert\leq l-\varepsilon
\} = \bigl\{ x\in D; \rho(x) \geq\varepsilon\bigr\}.
\end{equation}
%
%
%
\rule{\textwidth}{0.5pt}\hypertarget{algA1}\\
\textbf{Algorithm (A1) for} $\bolds{\delta=2}$. Let us fix a parameter
$0<\gamma<1$.\\
\textbf{Initialization:} Set $X(0)=(X_1(0),X_2(0))=0$, $\theta_0=0$,
$\Theta_0=0$, $A_0=\gamma^2l^2 e/2$.\\
\textbf{First step:} Let $(U_1,V_1,W_1)$ be a vector of three independent
random variables uniformly distributed on $[0,1]$. We set
\[
\cases{
\displaystyle\theta_1=A_0 U_1V_1,\qquad \Theta_1=\Theta_0+\theta_1,\vspace*{2pt}\cr
\displaystyle X(1)^\intercal=\bigl(X_1(1),X_2(1)\bigr)^\intercal=X(0)^\intercal+\psi
_{A_{0}}(\theta_1)\pmatrix{
\cos(2\pi W_1)\vspace*{2pt}\cr
\sin(2\pi W_1)},}
\]
where
%
\begin{equation}
\label{defdepsi}
\psi_a(t)=\sqrt{2t\log\frac{a}{t}}, \qquad  t\le a, a>0.
\end{equation}
At the end of this step we set $A_1=\gamma^2\rho(X(1))^2 e/2$.\\
\textbf{The $\bolds{n}$th step:} While $X (n-1)\in
D^\varepsilon$, simulate $(U_n,V_n,W_n)$ a vector of three independent
random variables uniformly distributed on $[0,1]$ and define
%
\begin{equation}\label{eqalgo-n}
\cases{\displaystyle
\theta_n=A_{n-1} U_nV_n,\qquad
\Theta_n=\Theta_{n-1}+\theta_n,\vspace*{2pt}\cr
\displaystyle X(n)^\intercal=\bigl(X_1(n),X_2(n)\bigr)^\intercal=X(n-1)^\intercal+\psi
_{A_{n-1}}(\theta_n)\pmatrix{\cos(2\pi W_n)\vspace*{2pt}\cr
\sin(2\pi W_n)}.}\hspace*{-35pt}
\end{equation}
At the end of this step we set $A_n=\gamma^2\rho(X(n))^2 e/2$.\\
When $X(n-1)\notin D^\varepsilon$ the algorithm is stopped: we set
$\theta_n=0$, $\Theta_n=\Theta_{n-1}$ and $X(n)=X(n-1)$.\\
\textbf{Outcome:} The hitting time $\Theta_n$ and the
exit position $X(n)$.\\
\rule{\textwidth}{0.5pt}
%
%

\begin{rem}
 The WoMS algorithm describes a $D$-valued Markov chain
$(X(n),  n\ge0)$. Each step corresponds to an exit problem for the
$2$-dimensional Brownian motion. If $X(n)=x$, then we focus our
attention to the exit problem of the ball centered in $x$ and of radius
$(\psi_a(t),  t\ge0)$: the exit location corresponds to $X(n+1)$ and
the exit time to $\theta_{n+1}$. Of course the choice of the parameter
$a$ plays a crucial role since the moving sphere has to belong to the
domain $D$ as time elapses. When the Markov chain $X$ is close to the
boundary $\partial D$, we stop the algorithm and obtain therefore a
good approximation of the exit problem of $D$.

Comparison with the classical ($\mathit{WoS}$) algorithm: at each step, the
$n$th step of the classical walk on spheres ($\mathit{WoS}$) is based on the
exit location and exit time, which are mutually independent, for the
Brownian paths exiting from a ball centered in $X(n-1)$ and with radius
$\gamma\rho(X(n-1))$. The exit location is uniformly distributed on
the sphere while the exit time is characterized by its Laplace
transform. Therefore, if one knows $X(n-1)$, then the diameter of the
sphere is deterministic. For the $\mathit{WoMS}$ the center of the ball should
also be $X(n-1)$, but the radius is random, smaller than $\gamma\rho
(X(n-1))$. The exit location will also be uniformly distributed on the
sphere, but the exit time will be much easier to simulate: in
particular, you do not need to evaluate the Bessel functions.
\end{rem}
The stochastic process $(X(n), n\ge0)$ is a homogeneous Markov chain
stopped at the first time it exits from $D^\varepsilon$. In the
following, we shall denote $N^\varepsilon$ this stopping time which
represents in fact the number of steps in the algorithm:
\[
N^\varepsilon=\inf\bigl\{n\ge0; X(n)\notin D^\varepsilon\bigr\}.
\]
We just notice that $X(N^\varepsilon)\notin D^\varepsilon$ by its
definition.\eject

Algorithm (\hyperlink{algA1}{A1}) is presented in the $2$-dimensional case. Of course we
can construct a generalization of this
procedure for the $\delta$-dimensional Bessel process when $\delta\in
\mathbb{N}^*$. For notational simplicity,
we use a slightly different method: instead of dealing with a Markov
chain $(X(n),  n\in\mathbb{N})$ living in $\mathbb{R}^\delta$ we shall
consider its squared norm,
which is also (surprisingly) a Markov chain. At each step, we shall
construct a couple of random variables $(\xi_n, \chi(n))$ associated
to an exit problem, the first coordinate corresponds to an exit time
and the second one to the norm of the exit location.

We introduce some notation: $\mathcal{S}^\delta$ represents the unit
ball in $\mathbb{R}^\delta$ and
$\pi_1\dvtx \mathbb{R}^\delta\to\mathbb{R}$ the projection on the first
coordinate.

\noindent\rule{\textwidth}{0.5pt}\hypertarget{algA2}\\
\textbf{Algorithm (A2).} Let us fix a parameter $0<\gamma<1$.\\
\textbf{Initialization:} Set $\chi(0)=0$, $\xi_0=0$, $\Xi_0=0$,
$A_0=(\gamma^2l^2 e/(\nu+1))^{\nu+1}\frac{\Gamma(\nu+1)}{2}$.\\
\textbf{The $\bolds{n}$th step:} While $\sqrt{\chi(n-1)}<l-\varepsilon$, we choose
$U_n$ a uniform distributed random vector on
$[0,1]^{\lfloor\nu\rfloor+2}$, $G_n$ a standard Gaussian random
variable and $V_n$ an uniformly distributed random vector on $\mathcal
{S}^\delta$. Consider
$U_n$, $G_n$ and $V_n$ independent. We set
%
\begin{eqnarray}\label{eqdefalgo}
\cases{
\displaystyle\xi_n=\biggl(\frac{A_{n-1}}{\Gamma(\nu+1)2^\nu}  U_n(1)\cdots U_n\bigl(\lfloor
\nu\rfloor+2\bigr)\biggr)^{{1}/{(\nu+1)}}\exp\biggl\{-\frac{\nu-\lfloor\nu
\rfloor}{\nu+1}G_n^2\biggr\}, \vspace*{2pt}\cr
\qquad\Xi_n=\Xi_{n-1}+\xi_n,\vspace*{2pt}\cr
\displaystyle\chi(n)=\chi(n-1)+2\pi_1(V_n)\sqrt{\chi(n-1)}\psi_{A_{n-1}}(\xi_n)+\psi
^2_{A_{n-1}}(\xi_n),}\hspace*{-35pt}
\end{eqnarray}
where
%
\begin{eqnarray}
\label{defdepsi-new}
\psi_a(t)&=&\sqrt{2t\log\frac{a}{\Gamma(\nu+1)t^{\nu+1}2^\nu}},
\nonumber
\\[-8pt]
\\[-8pt]
\nonumber
t&\le&
t_{\mathrm{max}}(a):=\left[\frac{a}{\Gamma(\nu+1)2^\nu}\right]^{{1}/{(\nu
+1)}}, \qquad a>0.
\end{eqnarray}
At the end of this step we set
\[
A_n=\bigl(\gamma^2\bigl(l-\sqrt{\chi(n)}\bigr)^2 e/(\nu+1)\bigr)^{\nu+1}\frac{\Gamma
(\nu+1)}{2}.
\]
When $\sqrt{\chi(n)}\ge l-\varepsilon$ the algorithm is stopped: we
then set $\xi_n=0$, $\Xi_n=\Xi_{n-1}$ and $\chi(n)=\chi(n-1)$.\\
\textbf{Outcome:} The hitting time $\Xi_n$ and the value of the Markov
chain $\chi(n)$.\\
\rule{\textwidth}{0.5pt}

It is obvious that for the particular dimension $\delta=2$, that is,
$\nu=0$, the stopping times obtained by Algorithms (\hyperlink{algA1}{A1}) and (\hyperlink{algA2}{A2})
have the same distribution.
Moreover, for each $n$, $\chi(n)$ has the same distribution as $\Vert
X(n) \Vert^2$. In other words, if the number of steps of (\hyperlink{algA1}{A1}) and (\hyperlink{algA2}{A2})
are identical in law, the number of steps will be denoted in both cases
$N^\varepsilon$.
%
%
%
\begin{thmm}
\label{thmalgo}
Set $\delta\in\mathbb{N}^*$. The number of steps $N^\varepsilon$ of
the Algorithm $\mathit{WoMS}$ \textup{(\hyperlink{algA2}{A2})} is almost surely finite. Moreover, there
exist constants $C_\delta>0$ and $\varepsilon_0(\delta)>0$, such that
\[
\mathbb{E}\bigl[N^\varepsilon\bigr]\le C_\delta|\log\varepsilon|
\qquad\mbox{for all } \varepsilon\le\varepsilon_0(\delta).
\]
\end{thmm}
%
%
%
%
\begin{thmm} \label{thmalgo-conv} Set $\delta\in\mathbb{N}^*$.
As $\varepsilon$ goes to zero, $\Xi_{N^\varepsilon}$ converges in
distribution toward
$\tau_l$, the hitting time of the $\delta$-dimensional Bessel process
(with cumulative distribution function $F$), which is almost surely
finite. Moreover, for any $\alpha>0$ small enough,
%
\begin{equation}
\label{eqthmencadr} \biggl(1-\frac{\eps}{\sqrt{2\alpha\pi}} \biggr)F^\eps(t-\alpha)
\le F(t)\le F^\eps(t)\qquad\mbox{for all } t>0,
\end{equation}
where $F^\eps(t):=\mathbb{P}(\Xi_{N^\varepsilon}\le t)$.
\end{thmm}
%
%
%
These results and the key ideas of the proofs are adapted from the
classical random walk on spheres ($\mathit{WoS}$); see \cite{Milstein-97}.

\begin{pf*}{Proof of Theorem \ref{thmalgo}}
\emph{Step} 1. Let us estimate the number of steps. Since $(\chi(n),
n\ge0)$ is a homogeneous Markov chain,
we introduce the operator $P_xf$ defined, for any nonnegative function
$f\dvtx \mathbb{R}_+\to\mathbb{R}_+$,
by
\[
P_xf:=\int_{\mathbb{R}_+}f(y)\mathbb{P}(x, \dint y),
\]
where $\mathbb{P}(x,\dint y)$ is the transition probability of the
Markov chain.
By definition, $\chi(n+1)$ depends only on $\chi(n)$, $V_n$ and $\xi
_n$. Let us note that,
by construction, $V_n$ and $\xi_n$ are independent. Moreover using the
result developed in the \hyperref[app]{Appendix}, the density of $\xi_n (\frac{2^\nu
\Gamma(\nu+1)}{A_{n-1}} )^{{1}/{(\nu+1)}}$ is given by
%
\begin{equation}
\label{eqdens} \mu(r)=\frac{(\nu+1)^{\nu+2}}{\Gamma(\nu+2)} r^\nu (-\log r
)^{\nu+1}\ind{[0,1]}(r).
\end{equation}
If we denote $\sigma^d$, the uniform surface measure on the unit sphere
in $\mathbb{R}^d$, we get
%
\begin{equation}
\label{eqpxf} \qquad P_xf=\int_{0}^1\int
_{\mathcal{S}^\delta}f \bigl(x+2\pi_1(u)\sqrt
{x}K(x,r)+K^2(x,r) \bigr)\mu(r)\,\dint r\, \sigma^\delta(\dint
u),
\end{equation}
with $K(x,r)$ defined by
%
\begin{equation}
\label{eqK} K(x,r)=\psi_{A} \biggl( \biggl[\frac{A}{2^\nu\Gamma(\nu+1)}
\biggr]^{
{1}/{(\nu+1)}}r \biggr),
\end{equation}
and $A$ depending on $x$ in the following way:
\[
A= \biggl(\frac{\gamma^2(l-\sqrt{x})^2 e}{\nu+1} \biggr)^{\nu+1}\frac{\Gamma
(\nu+1)}{2}.
\]
We can observe the following scaling property: $\psi_A(A^{{1}/{(\nu
+1)}}t)=A^{{1}/{(2\nu+2)}}\psi_1(t)$. Therefore the definition of $\psi
_1$ leads to
%
\begin{equation}
\label{eqdefK}\quad  K(x,r)=\gamma(l-\sqrt{x})\sqrt{\frac{er}{\nu+1}\log
\frac{1}{r^{\nu
+1}}}=\gamma(l-\sqrt{x})\sqrt{er(-\log r)}.
\end{equation}
\emph{Step} 2. Using classical potential theory for
discrete time Markov chains
(see, e.g., Theorem 4.2.3 in \cite{Norris-97}), we know that
\[
\phi(x)=\mathbb{E}_x \Biggl(\sum_{n=0}^{N^\varepsilon-1}g
\bigl(\chi(n)\bigr) \Biggr)
\]
satisfies, for any nonnegative function $g$,
%
\begin{equation}
\label{eqequat} \cases{ %
 \phi(x)=P_x
\phi+g(x),&\quad $0\le x < (l-\varepsilon)^2,$
\vspace*{2pt}\cr
\phi\bigl((l- \varepsilon)^2\bigr)=0.}
\end{equation}
In particular, for $g=1$, we obtain that $\phi(x)=\mathbb
{E}_x[N^\varepsilon]$.
In order to get an upper-bound for the averaged number of steps, it
suffices to apply
a comparison result.
Let us first define the constant $C_\delta$,
%
\begin{equation}
\label{Cdelta} C_\delta= \biggl(\frac{\nu+1}{\nu+2} \biggr)^{\nu+2}
\frac{e}{\Gamma
(\nu+2)}\frac{1}{2\delta}\sigma^{\delta}\bigl(S^\delta
\bigr).
\end{equation}
We choose the function
%
\begin{equation}\qquad
U^\varepsilon(x)=\bigl\{\log\bigl((l-\sqrt{x})/\varepsilon\bigr)-\log(1-
\gamma)\bigr\} /\bigl(C_\delta\gamma^2\bigr), \qquad 0\leq x <
l^2,
\end{equation}
which satisfies $U^\varepsilon(x)\ge P_xU^\varepsilon+1$, for all
$0<x<(l-\varepsilon)^2$ (see Lemma \ref{lemineg-log} for the
definition of the constant and for the inequality) and $U^\varepsilon
(x)\ge0$ for all $0<x<(l-\varepsilon)^2$. A classical comparison result
related to the potential theory (see, e.g., Theorem 4.2.3 in
\cite{Norris-97}) implies that $\mathbb{E}_x[N^\varepsilon]\le
U^\varepsilon(x)$ for all $x\in[0,(l-\varepsilon)^2]$ and consequently
leads to the announced statement.
\end{pf*}
%
%
%
\begin{lemma}
\label{lemineg-log} Let us define, for small $\varepsilon>0$,
$U^\varepsilon(x)=\{\log((l-\sqrt{x})/\varepsilon)-\log(1-\gamma)\}
/(C_\delta\gamma^2)$ for $x\in[0,l^2[$ and where the constant
$C_\delta$ is given by (\ref{Cdelta}) and $\gamma$ is related to the
definition of the $\mathit{WoMS}$. Then, for any $x\in]0,(l-\varepsilon)^2[$,
the following inequality yields
\[
P_xU^\varepsilon-U^\varepsilon(x)\le-1.
\]
We recall that $P_xU^\varepsilon$ is defined by \eqref{eqpxf} and
\eqref{eqdefK}.
\end{lemma}

\begin{pf} We will split the proof into several steps.

\emph{Step} 1. First of all, we observe that $U^\varepsilon\ge-\log
(1-\gamma)/(C_\delta\gamma^2)$ in the domain $[0,(l-\varepsilon)^2]$.
Let us consider now $\chi(0)=x\in[0,(l-\varepsilon)^2]$ and $y$ in the
support of the law of $\chi(1)$ and let us prove that $U^\varepsilon
(y)\ge0$. By the definition of $\chi(1)$ we obtain
\[
\chi(1)\le\sup_{y\in[-1,1], t\in[0,t_\mathrm{max}(A)]} \bigl(x+2y\sqrt{x}\psi
_{A}(t)+\psi^2_A(t)\bigr),
\]
where $A= (\gamma^2(l-\sqrt{x})^2 e/(\nu+1) )^{\nu+1}\frac{\Gamma
(\nu+1)}{2}$ and both $\psi_A$ and $t_\mathrm{max}$ are defined by~\eqref
{defdepsi-new}. The right-hand side of the preceding inequality is
increasing with respect to $y$ so that
\[
\chi(1)\le \Bigl(\sqrt{x}+\sup_{t\in[0,t_\mathrm{max}(A)]} \psi_A(t)
\Bigr)^2.
\]
Furthermore, for $a>0$ the maximum of the function $\psi_a$ is reached
for $t_\mathrm{max}(a)=\frac{1}{e}   ( \frac{a}{\Gamma(\nu+1)2^\nu}
)^{{1}/{(\nu+1)}}$ and is equal to
%
\begin{equation}
\label{eqmaxcalcul} \sup_{t\in[0,t_\mathrm{max}(a)]} \psi_a(t)= \biggl
\{\frac{2(\nu+1)}{e} \biggl(\frac{a}{\Gamma(\nu+1)2^\nu} \biggr)^{{1}/{(\nu+1)}} \biggr
\}^{1/2}.
\end{equation}
Finally using the definition of $A$ and the inequality $x\le
(l-\varepsilon)^2$, we find the following lower bound:
\[
l-\sqrt{\chi(1)}\ge(l-\sqrt{x}) (1-\gamma)\ge\varepsilon(1-\gamma).
\]
We can therefore conclude that, for any $y$ in the support of the law
of $\chi(1)$ (even for $y\ge(l-\varepsilon)^2$), $U^\varepsilon(y)\ge0$
which ensures that $U^\varepsilon$ is well defined and nonnegative in
the domain of the operator $P_x$.

\emph{Step} 2. Furthermore the Taylor expansion yields
%
\begin{eqnarray}
\label{eqtaylor}\qquad U^\eps(y)\le U^\eps(x)+\frac{\sqrt{x}-\sqrt{y}}{C_\delta\gamma^2
(l-\sqrt{x})}-
\frac{(\sqrt{x}-\sqrt{y})^2}{2C_\delta\gamma^2(l-\sqrt {x})^2}+\frac{(\sqrt{x}-\sqrt{y})^3}{3C_\delta\gamma^2(l-\sqrt {x})^3},
\nonumber
\\[-8pt]
\\[-8pt]
\eqntext{x, y\in[0,l^2[.}
\end{eqnarray}
If $\chi(0)=x$ and $y$ is in the support of the random variable $\chi
(1)$, then
\begin{eqnarray*}
\sqrt{y}-\sqrt{x}&=&\sqrt{x+2\pi_1(u)\sqrt{x}K(x,r)+K^2(x,r)}-
\sqrt{x}
\\
&\ge&\pi_1(u)K(x,r).
\end{eqnarray*}
By expansion \eqref{eqtaylor} and the definition of the operator $P_x$
given by \eqref{eqpxf}, the following upper-bound for the operator
$P_x$ holds:
\begin{eqnarray*}
P_xU^\eps&=&\int_{0}^1
\int_{\mathcal{S}^\delta}U^\varepsilon \bigl(x+2\pi _1(u)
\sqrt{x}K(x,r)+K^2(x,r) \bigr)\mu(r)\,\dint r \,\sigma^\delta(
\dint u),
\\
&\le& U^\eps(x)-\int_{0}^1\int
_{\mathcal{S}^\delta}\frac{\pi
_1(u)K(x,r)}{C_\delta\gamma^2(l-\sqrt{x})} \mu(r)\,\dint r \,\sigma^\delta (
\dint u)
\\
&&{} - \int_{0}^1\int_{S^\delta_+}
\frac{\pi
_1^2(u)K^2(x,r)}{2C_\delta\gamma^2(l-\sqrt{x})^2} \mu(r)\,\dint r \,\sigma ^\delta(\dint u)
\\
&&{} -\int_{0}^1\int_{\mathcal{S}^\delta}
\frac{\pi
_1^3(u)K^3(x,r)}{3C_\delta\gamma^2(l-\sqrt{x})^3} \mu(r)\,\dint r \,\sigma ^\delta(\dint u),
\end{eqnarray*}
where
%
\begin{equation}
\label{Sdelta+} S^\delta_+: =\bigl\{u\in\mathcal{S}^\delta
\dvtx \pi_1 (u) > 0\bigr\}.
\end{equation}
Due to symmetry properties, the first and the third integral terms
vanish. Then~\eqref{eqdefK} leads to
\[
P_xU^\eps\le U^\eps(x)- \frac{I}{C_\delta}
\int_{S^\delta_+} \pi _1^2(u)
\sigma^\delta(\dint u)
\]
with
\[
I=\frac{(\nu+1)^{\nu+2}e}{2\Gamma(\nu+2)}\int_0^1r^{\nu+1}(-
\log r)^{\nu+2} \,\dint r.
\]
The description of the probability density function in the \hyperref[app]{Appendix}
leads to the following explicit value:
\[
I= \biggl( \frac{\nu+1}{\nu+2} \biggr)^{\nu+2}\frac{e}{\Gamma(\nu+2)}.
\]
In order to complete the proof, it suffices to choose the particular
constant given by (\ref{Cdelta}) after noticing that
\[
\int_{S^\delta_+} \pi_1^2(u)
\sigma^\delta(\dint u)= \frac
{1}{2\delta}\sigma^\delta
\bigl(S^\delta\bigr).
\]
\upqed\end{pf}
%
%
%
\begin{pf*}{Proof of Theorem \ref{thmalgo-conv}}
The proof is split in two parts. First, the steps of the algorithm
and the hitting time of the Bessel process of index $\nu$ shall be
related to stopping times of a $\delta$-dimensional Brownian motion
($\nu=\frac{\delta}{2}-1$). Second, we point out that the corresponding
stopping times are close together by evaluating deviations of the
Brownian paths.

\emph{Step} 1. Let $\mathbf{B}=(B^{(1)},B^{(2)},\ldots,B^{(\delta)})$
be a $\delta$-dimensional Brownian motion. Then the norm of $\mathbf{B}$
has the same distribution as a Bessel process of index $\nu$; see, for
instance, \cite{Revuz-Yor-99}. Hence the first hitting time $\tau_l$ is
identical in law to the stopping time
\[
\mathbb{T}_l=\inf\{t\ge0; \mathbf{B}_t\notin D \},
\]
where $D$ is the Euclidean ball centered at the origin and of radius
$l$. We introduce then a procedure in order to come close to $\mathbb
{T}_l$. For the first step we shall focus our attention to the first
exit time of a moving sphere centered at the origin and of radius $\psi
_a(t)$ defined by \eqref{defdepsi-new}, we denote $\hat{\xi}_1$ this
stopping time. Of course this moving sphere should always stay in $D$,
so we choose $a$ such that the maximum of $\psi_a$ stays smaller than
$l$. By \eqref{eqmaxcalcul}, we get
\[
\sup_{t\le a}\psi_a(t)< l\quad\Longleftrightarrow\quad a<
\frac{\Gamma(\nu
+1)}{2} \biggl(\frac{el^2}{\nu+1} \biggr)^{\nu+1}.
\]
For $a=A_0=\frac{\Gamma(\nu+1)}{2} (\frac{e\gamma^2l^2}{\nu+1}
)^{\nu+1}$ with a parameter $\gamma<1$, the condition is satisfied,
$\sup_{t\le a}\psi_a(t)=\gamma^{2\nu+2} l<l$. Let us describe the law
of $(\hat{\xi}_1,\mathbf{B}_{\hat{\xi}_1})$. The norm of the Brownian
motion is identical in law with the Bessel process; therefore
Proposition~\ref{propfirstboundary} implies that the density function
of $\hat{\xi}_1$ is given by
\eqref{densitetau} with $a$ replaced by $A_0$. Using the law described
in Proposition \ref{propappend}, we can prove that $\hat{\xi}_1$ has
the same distribution as
\[
\biggl(\frac{ A_0}{\Gamma(\nu+1)2^\nu} \biggr)^{{1}/{(\nu+1)}} e^{-Z},
\]
where $Z$ is Gamma distributed with parameters\vspace*{-1.5pt} $\alpha=\nu+2$ and $\beta
=\frac{1}{\nu+1}$.
By construction we deduce that $\hat{\xi}_1 \stackrel{(d)}{=}
\xi_1$ where $\xi_1$ is defined in the Algorithm $\mathit{WoMS}$~(\hyperlink{algA2}{A2}). Knowing
the stopping time, we can easily describe the exit location since the
Brownian motion is rotationnaly invariant: $\mathbf{B}_{\hat{\xi}_1}$ is
then uniformly distributed on the sphere of radius $\psi_{A_0}(\hat{\xi
}_1)$. Hence
\[
\bigl(\hat{\xi}_1,\Vert\mathbf{B}_{\hat{\xi}_1}\Vert\bigr)\stackrel{\mathrm{(d)}} {=}
\bigl(\xi _1,\chi(1)\bigr)\quad \mbox{and}\quad \hat{\xi}_1<
\mathbb{T}_l.
\]
By this procedure we can construct a sequence of stopping times $(\hat
{\xi}_n,  n\ge1)$ and define $\hat{\Xi}_n=\hat{\xi}_1+\cdots+\hat{\xi
}_n$; $\hat{\Xi}_n$ is the first time after $\hat{\Xi}_{n-1}$ such that
the Brownian motion exits from a sphere centered in $\mathbf{B}_{\hat{\Xi
}_{n-1}}$ of radius $\psi_{a_n}$ initialized at time $\hat{\Xi}_{n-1}$.
See Figure \ref{fig1}. The moving sphere should stay in the domain $D$, so we choose
\[
a_n= \bigl(\gamma^2 (l-\sqrt{\mathbf{B}_{\hat{\Xi}_{n-1}}}
)^2 e/(\nu +1) \bigr)^{\nu+1}\frac{\Gamma(\nu+1)}{2}.
\]

\begin{figure}

\includegraphics{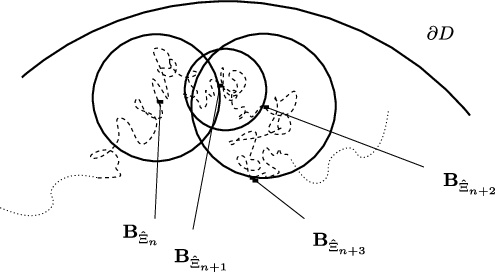}

\caption{Walk on moving spheres.}\label{fig1}\vspace*{-3pt}
\end{figure}

%
%

Using the same arguments as before and by the Markov property for the
Brownian motion, we obtain the identities in law
\[
(a_n, n\ge1)\stackrel{\mathrm{(d)}} {=} (A_n, n\ge1),\qquad \bigl(\hat{\Xi
}_n, \Vert\mathbf{B}_{\hat{\Xi}_n}\Vert \bigr)_{n\ge1}
\stackrel{\mathrm{(d)}} {=} \bigl(\Xi_n, \chi(n) \bigr)_{n\ge1}
\]
with $\hat{\Xi}_n<\mathbb{T}_l$ and $\Xi_n$, $A_n$, $\chi(n)$ defined
in the Algorithm $\mathit{WoMS}$ (\hyperlink{algA2}{A2}). Consequently defining $\hat{N}^\varepsilon
=\inf\{ n\ge0; \mathbf{B}_{\hat{\Xi}_n}\notin D^\eps\}$, the following
identity yields
%
\begin{equation}
\label{eqidenlaw}\bigl ( \hat{\Xi}_{\hat{N}^\eps}, \Vert\mathbf{B}_{\hat{\Xi}_{\hat{N}^\eps}}
\Vert \bigr)\stackrel{\mathrm{(d)}} {=} \bigl( \Xi_{N^\eps}, \chi\bigl(N^\eps
\bigr) \bigr) \quad\mbox{and}\quad \hat{\Xi}_{\hat{N}^\eps}<\mathbb{T}_l.
\end{equation}
\emph{Step} 2. Let us now estimate the difference between $\hat{\Xi
}_{\hat{N}^\eps}$ and $\mathbb{T}_l$. By \eqref{eqidenlaw} we first deduce
%
\begin{equation}
\label{equpper-bound} F(t):=\mathbb{P}(\tau_l\le t)=\mathbb{P}(
\mathbb{T}_l\le t)\le F^\eps (t):=\mathbb{P}(
\Xi_{N^\eps}\le t),\qquad t>0.
\end{equation}
Furthermore, for any small $\alpha>0$,
%
\begin{eqnarray}
\label{eqdecomp-under} 1-F(t)&=&\mathbb{P}( \mathbb{T}_l> t, \hat{
\Xi}_{\hat{N}^\eps}\le t-\alpha)+\mathbb{P}( \mathbb{T}_l> t, \hat{
\Xi}_{\hat{N}^\eps}> t-\alpha)
\nonumber
\\
&\le&\mathbb{P}( \mathbb{T}_l> t, \hat{\Xi}_{\hat{N}^\eps}\le t-
\alpha)+\mathbb{P}(\hat{\Xi}_{\hat{N}^\eps}> t-\alpha)
\\
&\le&\mathbb{P}( \mathbb{T}_l> t, \hat{\Xi}_{\hat{N}^\eps}\le t-
\alpha)+1-F^\eps(t-\alpha).\nonumber
\end{eqnarray}
At time $\hat{\Xi}_{\hat{N}^\eps}$ the Brownian motion is in the $\eps
$-neighborhood of the boundary $\partial D$, hence
$l - \Vert\mathbf{B}_{\hat{\Xi}_{\hat{N}^\eps}}\Vert \le \eps$. Using the
strong Markov property, we obtain
%
\begin{equation}
\label{eqinterm} \mathbb{P}( \mathbb{T}_l> t, \hat{
\Xi}_{\hat{N}^\eps}\le t-\alpha )\le F^\eps(t-\alpha)\sup
_{y\in D\setminus D^\eps}\mathbb{P}_y(\mathbb {T}_l>
\alpha).
\end{equation}
Since the Brownian motion is rotationally invariant, it suffices to
choose $y=(l-\varepsilon,0,\ldots,0)$. Due to the convexity of $D$, the
following upper-bound holds:
%
\begin{equation}
\label{eqlowboununif} \mathbb{P}_y(\mathbb{T}_l>\alpha)
\le\mathbb{P}_0\Bigl(\sup_{0\le t\le\alpha
}
\overline{B}_t^{(1)}<\eps\Bigr)=\mathbb{P}_0
\bigl(2\bigl\vert\overline{B}_\alpha ^{(1)}\bigr \vert<\eps\bigr)\le
\frac{\eps}{\sqrt{2\alpha\pi}}.
\end{equation}
Combining \eqref{equpper-bound} for the upper-bound and \eqref
{eqdecomp-under}, \eqref{eqinterm} and \eqref{eqlowboununif} for the
lower-bound yields the announced estimation \eqref{eqthmencadr}.
\end{pf*}
%
%
\subsubsection{\texorpdfstring{The first time the Bessel process of index $\nu$ hits a
decreasing curved boundary}{The first time the Bessel process of index nu hits a
decreasing curved boundary}}
The algorithm developed in the previous paragraph can be\break adapted to the
problem of hitting a deacreasing curved boundary. Let us define
%
\begin{equation}
\label{deftau-line}\qquad \tau=\inf\bigl\{t\ge0\dvtx Z_t^{\delta,0}=l(t)
\bigr\}\qquad \mbox{where } l \mbox{ is decreasing and } l(0)>0.
\end{equation}
%
\begin{assu}\label{assu}
There exists a constant $\Delta_{{\mathrm{min}}}>0$ which bounds the
derivative of $l$
\[
l'(t)\ge-\Delta_\mathrm{min} \qquad\forall t\ge0.
\]
\end{assu}
The procedure then also consists in building a $\mathit{WoMS}$ which reaches a
neighborhood of the boundary. But instead of dealing with a fixed
boundary as in Section~\ref{sect0-horiz}, that is a ball of radius
$l$, we shall in this section introduce the following moving boundary:
the ball centered in the origin and of radius $l(t)$. The arguments
developed in order to prove Theorems \ref{thmalgo} and \ref
{thmalgo-conv} will be adapted to this new context.
%
%

\noindent\rule{\textwidth}{0.5pt}\hypertarget{algA3}\\
\textbf{Algorithm (A3):}\\
Let us define the following positive constants:
%
\begin{equation}
\label{constante}
L=\max\bigl(l(0),\Delta_{\mathrm{min}}, \sqrt{\nu+1}\bigr),\qquad \kappa=\frac
{2^\nu}{5^{\nu+1}L^{2\nu+2}} \Gamma(\nu+1).
\end{equation}
\textbf{Initialization:} Set $\chi(0)=0$, $\xi_0=0$, $\Xi_0=0$, $A_0=\kappa
( l(0)- \sqrt{\chi(0)})^{2(\nu+1)}$.\\
\mathversion{bold}
\textbf{The $n$th step:}\mathversion{normal} While the condition
\[
l(\Xi_{n-1})-\sqrt{\chi(n-1)}>\varepsilon
\]
[denoted by $\mathbb{C}(n-1)$] holds, we simulate $U_n$ an uniform
distributed random vector on
$[0,1]^{\lfloor\nu\rfloor+2}$, $G_n$ a standard Gaussian random
variable and $V_n$ a uniformly distributed random vector on $\mathcal
{S}^\delta$.
$U_n$, $G_n$ and $V_n$ have to be independent. We then construct $(\xi
_n,\Xi_n,\chi(n))$ using \eqref{eqdefalgo}. At the end of this step we
set $A_n=\kappa( l(\Xi_n)-\sqrt{\chi(n)})^{2(\nu+1)}$.\\
The algorithm stops when $\mathbb{C}(n-1)$ is not longer satisfied: we
set $\xi_n=0$ and so $\Xi_n=\Xi_{n-1}$ and $\chi(n)=\chi(n-1)$.\\
\textbf{Outcome} The exit position $\chi(n)$ and the exit time.\\
\rule{\textwidth}{0.5pt}

Let us note that the stochastic process $(\chi(n), n\ge0)$ is not a
Markov chain since the sequence $(A_n)_{n\ge0}$ depends on both $(\Xi
_n, \chi(n))$. That is why we define the following Markov chain:
\[
R_n:=\bigl(\Xi_n,\chi(n)\bigr)\in\mathbb{R}_+^2
\]
stopped at the first time the condition $\mathbb{C}(n)$ is not
satisfied. In the following, we shall denote $N^\varepsilon$ this
stopping time (number of steps of the algorithm):
\[
N^\varepsilon=\inf\bigl\{n\ge0; l(\Xi_n)-\sqrt{\chi(n)}\le\varepsilon\bigr\}.
\]

%
%
%
\begin{thmm}
\label{thmalgoline}
The number of steps $N^\varepsilon$ of the Algorithm $\mathit{WoMS}$ \textup{(\hyperlink{algA3}{A3})} is
almost surely finite. Moreover, there exist a constant $C_\delta>0$ and
$\varepsilon_0(\delta)>0$, such that
\[
\mathbb{E}\bigl[N^\varepsilon\bigr]\le C_\delta|\log\varepsilon|\qquad
\mbox{for all } \varepsilon\le\varepsilon_0(\delta).
\]
\end{thmm}
%
%
%
%
\begin{thmm} \label{thmalgo-convline}
As $\varepsilon$ goes to zero, $\Xi_{N^\varepsilon}$ converges in
distribution toward $\tau$ defined by \eqref{deftau-line} (with
cumulative distribution function $F$), which is almost surely finite.
Moreover, for any $\alpha>0$ small enough,
%
\begin{equation}
\label{eqthmencadrline} \biggl(1-\frac{\eps}{\sqrt{2\alpha\pi}} \biggr)F^\eps(t-
\alpha)\le F(t)\le F^\eps(t)\qquad\mbox{for all } t>0,
\end{equation}
where $F^\eps(t):=\mathbb{P}(\Xi_{N^\eps}\le t)$.
\end{thmm}
\begin{pf*}{Proof of Theorem \ref{thmalgoline}}
The proof is based mainly on arguments already presented in Theorem \ref
{thmalgo}. So we let the details of the proof to the reader and focus
our attention to the main ideas.

(1) The process $(\Xi_n, \chi(n))$ is a homogeneous Markov chain
and the associated operator is given by
%
\begin{equation}
\label{eqnew-operator} P_{t,x}f:=\int_{(s,y)\in\mathbb{R}_+^2}f(s,y)
\mathbb{P} \bigl((t,x),(\dint s,\dint y) \bigr),
\end{equation}
where $f$ is a nonnegative function and $\mathbb{P} ((t,x),(\dint
s,\dint y) )$ is the transition probability of the chain. The chain
starts with $(\Xi_0,\chi(0))=(0,0)$ and is stopped the first time when
$l(\Xi_n)-\sqrt{\chi(n)}\le\varepsilon$. Classical potential theory
ensures that
\[
\phi(t,x)=\mathbb{E}_{t,x} \Biggl( \sum_{n=0}^{N^\varepsilon-1}g
\bigl(\Xi _n,\chi(n)\bigr) \Biggr)
\]
is solution of the following equation:
%
\begin{eqnarray}
\label{eqsysline} \cases{ %
 \phi(t,x)=P_{t,x}
\phi+g(t,x),&\quad $(t,x)\in D^\varepsilon,$
\vspace*{2pt}\cr
\phi(t,x)=0, &\quad$\forall(t,x)\in\partial D^\varepsilon,$ }
\end{eqnarray}
with $D^\varepsilon=\{(t,x)\in\mathbb{R}_+^2\dvtx  l(t)-\sqrt{x}\le
\varepsilon\}$.
For the particular choice $g=1$, we obtain $\phi(t,x)=\mathbb
{E}_{t,x}[N^\varepsilon]$, and
therefore the averaged number of step is given by $\phi(0,0)$.

(2) In order to point out an upper-bound for the averaged number
of steps,
we use a comparison result: we are looking for a function $U(t,x)$ such that
%
\begin{eqnarray}
\label{sysavec1} \cases{ %
 U(t,x)\ge P_{t,x}U+1,&\quad
$\forall(t,x)\in D^\varepsilon,$
\vspace*{2pt}\cr
U(t,x)\ge0, &\quad $\forall(t,x)\in\partial D^\varepsilon.$ }
\end{eqnarray}
For such a particular function, we can deduce $\phi(t,x)\le U(t,x)$.
Let us define
\[
U(t,x)=c\log \biggl( \frac{l(t)-\sqrt{ x}}{\varepsilon} \biggr)1_{\{
l(t)-\sqrt{x} \ge0 \}},
\]
with some constant $c>0$ which shall be specified later on. The
positivity assumption on the boundary $\partial D^ \varepsilon$
is trivial. Moreover since $l$ is a decreasing function, \eqref
{eqnew-operator} implies
%
\begin{eqnarray}
\label{eqineg-inte} P_{t,x}U&=&\int_{(s,y)\in\mathbb{R}_+^2}U(s,y)
\mathbb{P} \bigl((t,x),(\dint s,\dint y) \bigr)
\nonumber
\\[-8pt]
\\[-8pt]
\nonumber
&\le& \int_{(s,y)\in\mathbb{R}_+^2}U(t,y)
\mathbb{P} \bigl((t,x),(\dint s,\dint y) \bigr).
\end{eqnarray}
By using the Taylor expansion, we get
%
\begin{eqnarray}
\label{eqtaylor-} \qquad U(t,y)\le U(t,x)-c \frac{\sqrt{y}-\sqrt{x}}{l(t)-\sqrt{x}}-\frac
{c}{2}
\frac{(\sqrt{y}-\sqrt{x})^2}{(l(t)-\sqrt{x})^2}-\frac{c}{3}\frac
{(\sqrt{y}-\sqrt{x})^3}{(l(t)-\sqrt{x})^3},
\nonumber
\\[-8pt]
\\[-8pt]
\eqntext{(x,y)\in
\mathbb{R}_+^2.}
\end{eqnarray}
Using similar arguments and similar bounds as those presented in Lemma
\ref{lemineg-log}, the odd powers in the Taylor expansion do not play
any role in the integral \eqref{eqineg-inte}. Therefore we obtain
\begin{eqnarray*}
P_{t,x}U&\le& U(t,x)-\frac{c}{2} \int_{(s,y)\in\mathbb{R}_+^2}
\frac
{(\sqrt{y}-\sqrt{x})^2}{(l(t)-\sqrt{x})^2}\mathbb{P} \bigl((t,x),(\dint s,\dint y) \bigr)
\\
&\le& U(t,x)- \frac{c}{2}\int_{0}^1\int
_{S^\delta_+}\frac{\pi
_1^2(u)K^2(x,r)}{(l(t)-\sqrt{x})^2} \mu(r)\,\dint r \,\sigma^\delta(
\dint u),
\end{eqnarray*}
where $S^\delta_+$ is given in \eqref{Sdelta+}, and $K$ is defined by
\eqref{eqK} with $A=\kappa(l(s)-\sqrt{x})^{2(\nu+1)}$. We have now
\begin{eqnarray*}
P_{t,x}U&\le& U(t,x)-\frac{c(\nu+1)}{2} \biggl( \frac{2K}{\Gamma(\nu+1)}
\biggr)^{{1}/{(\nu+1)}} \biggl( \int_{\mathcal{S}^\delta_+}
\pi_1^2(u) \sigma ^\delta(\dint u) \biggr)\\
&&\hspace*{44pt}{}\times \biggl(
\int_0^1 r(-\log r)\mu(r)\,\dint r \biggr).
\end{eqnarray*}
An appropriate choice of the constant $c$ leads to \eqref{sysavec1}.
Finally we get
\[
\mathbb{E}\bigl[N^\varepsilon\bigr]\le U(0,0)= c \log\bigl(l(0)/\varepsilon
\bigr).
\]
\upqed\end{pf*}

\begin{pf*}{Proof of Theorem \ref{thmalgo-convline}}
The arguments are similar to those developed for Theorem \ref
{thmalgo-conv}, and the extension of the convergence result to curved
boundaries is straightforward. That is why we shall not repeat the
proof, but just focus our attention on the only point which is quite
different. We need to prove that the Markov chain $R_n=(\Xi_n,\chi(n))$
stays in the domain $D^0=\{(t,x)\dvtx  0\le x\le l^2(t) \}$ so that the
hitting time $\tau$ defined by \eqref{deftau-line} satisfies $\tau>\Xi
_{N^\varepsilon}$. In other words, if the Markov chain $R_n=(\Xi_n,\chi
(n))$ for the $n$th step is equal to $(s,x)$, then $R_{n+1}$ should
belong to $\{(t,x)\dvtx  t\ge s,  x\le l^2(t)\}$. In the $\mathit{WoMS}$ setting,
for $t\ge s$, this means that the ball centered in $x$ and of
time-dependent radius $\psi_A(t-s)$ always belongs as time elapses to
the ball centered in $0$ of radius $l(t)$. We recall that
\[
A=\kappa\bigl(l(s)-\sqrt{x}\bigr)^{2(\nu+1)}.
\]
Therefore we shall prove that
%
\begin{equation}
\label{eqamontrer} \forall t\ge s\qquad \psi_A(t-s)+\sqrt{x}\le l(t).
\end{equation}
In fact, due to Assumption \ref{assu} and the definition of $\psi_A$,
it suffices to obtain
%
\begin{equation}
\label{eprou} \psi_A(t-s)\le l(s)-\sqrt{x}-\Delta_{{\mathrm{min}}}(t-s)\qquad
\forall s\le t\le s+W^2,
\end{equation}
where
\begin{eqnarray*}
W&=& \biggl(\frac{A}{\Gamma(\nu+1)2^\nu} \biggr)^{{1}/{(2\nu+2)}}= \biggl(\frac{\kappa}{\Gamma(\nu+1)2^\nu}
\biggr)^{{1}/{(2\nu+2)}}\bigl(l(s)-\sqrt {x}\bigr)\\
&=&\frac{1}{L\sqrt{5}} \bigl(l(s)-
\sqrt{x}\bigr).
\end{eqnarray*}
Due to the definition of the constant $L$, we have
\begin{eqnarray*}
0&\le& W\le\frac{1}{2\Delta_{{\mathrm{min}}}} \frac{2(l(s)-\sqrt{x})\Delta
_\mathrm{min}}{\sqrt{({(2\nu+2)}/{e})+4(l(s)-\sqrt{x})\Delta_\mathrm{min}}}\\
&\le&\frac{1}{2\Delta_{{\mathrm{min}}}} \biggl\{
\sqrt{ \frac{2\nu
+2}{e}+4\bigl(l(s)-\sqrt{x}\bigr)\Delta_\mathrm{min}} -
\sqrt{ \frac{2\nu+2}{e} } \biggr\}.
\end{eqnarray*}
The right-hand side of the preceding inequality is the positive root of
the polynomial function $P(X)=\Delta_\mathrm{min} X^2+\sqrt{2(\nu+1)/e}
X-(l(s)-\sqrt{x})$. We deduce that $P(W)\le0$.
By \eqref{eqmaxcalcul} and $P(W)\le0$, we obtain
\begin{eqnarray*}
\sup_{t\ge s}\psi_A(t-s) &=& \biggl(
\frac{2(\nu+1)}{e} \biggr)^{1/2} W
\\
&\le& l(s)-\sqrt{x} -\Delta_\mathrm{min} W^2
\\
&\le &l(s)-\sqrt{x}-\Delta_{{\mathrm{min}}}(t-s)\qquad \forall s\le t\le
s+W^2.
\end{eqnarray*}
Finally we have proved \eqref{eprou} and so \eqref{eqamontrer}.
\end{pf*}

If Assumption \ref{assu} is not satisfied, then it is difficult to have
a general description of an iterated procedure in order to simulate
hitting times. However the particular form of the function $\psi_a$
defined by \eqref{defdepsi-new} permits us to describe a $\mathit{WoMS}$
algorithm for the square root boundaries.
Let us therefore consider the following functions:
%
\begin{equation}
\label{eqdefde2} \psi_a(t)=\sqrt{2t\log\frac{a}{\Gamma(\nu+1)t^{\nu+1}2^\nu}}\quad \mbox
{and}\quad f(t)=\sqrt{r-ut},
\end{equation}
well defined for $t\le t_0:=\min ( \alpha^{{1}/{(\nu+1)}},\frac
{r}{u}  )$ where $\alpha=a(\Gamma(\nu+1)2^\nu)^{-1}$.\eject

The algorithm is essentially based on the following result (the
constants $r$ and~$u$ associated with the hitting problem of a square
root boundary for the Bessel process shall be specified in the proof of
Proposition \ref{propcurved}).
%
\begin{lemma}
\label{lemcomparai}
Let us define
%
\begin{equation}
\label{eqdefdeF} F_\nu(r,u)=\frac{1}{2} \biggl(
\frac{er}{\nu+1} \biggr)^{\nu+1}\Gamma(\nu +1) e^{-u/2},\qquad r>0,
u>0.
\end{equation}
If $a=F_\nu(r,u)$, then
%
\begin{equation}
\label{comparai} \psi_a(t)\le f(t)\qquad \mbox{for all } 0\le t\le
\alpha^{{1}/{(\nu+1)}}.
\end{equation}
\end{lemma}

\begin{pf}
We are looking for a particular value $a$ depending on both $r$ and
$u$ such that the following bound holds:
$\psi_a(t)\le f(t)$, for all $0\le t\le t_0$.
Since $t\le t_0$, it suffices to prove that
\[
2t\log\frac{\alpha}{t^{\nu+1}}\le r-ut \quad\Longleftrightarrow\quad g(t):=t \biggl( 2\log
\frac{\alpha}{t^{\nu+1}}+u \biggr)\le r.
\]
Let us compute the maximum of the function $g$ on the interval
$[0,t_0]$, with $t_0$ fixed,
\[
g'(t)=2\log\frac{\alpha}{t^{\nu+1}}+u-2(\nu+1).
\]
We have
\[
g'(t)=0\quad \Longleftrightarrow\quad\log\frac{\alpha}{t^{\nu+1}}=\nu+1-
\frac
{u}{2}\quad\Longleftrightarrow\quad t^{\nu+1}=\alpha\exp \biggl\{
\frac{u}{2}-\nu-1 \biggr\}.
\]
In other words the maximum of the function $g$ is reached for
\[
t_{\mathrm{max}}=\alpha^{{1}/{(\nu+1)}}\exp \biggl\{ \frac{u}{2(\nu+1)}-1
\biggr\}
\]
and is equal to
\[
g(t_{\mathrm{max}})=g_{\max}=2(\nu+1)\alpha^{{1}/{(\nu+1)}}e^{({u}/{(2(\nu+1))})-1}.
\]
Choosing $g_{\max}\le r$ we obtain in particular \eqref{comparai},
which means
\[
\alpha\le \biggl(\frac{er}{2(\nu+1)} \biggr)^{\nu
+1}e^{-u/2}\quad
\Longleftrightarrow \quad a\le\frac{1}{2} \biggl(\frac{er}{\nu
+1}
\biggr)^{\nu+1}\Gamma(\nu+1)e^{-u/2}.
\]
For $a_0= \frac{1}{2} (\frac{er}{\nu+1} )^{\nu+1}\Gamma(\nu
+1)e^{-u/2}$, we get \eqref{comparai} since $t_0=\alpha^{{1}/{(\nu+1)}}$.
\end{pf}

The aim is now to construct an algorithm which permits us to
approximate the hitting time of the square root boundary. Therefore we
consider\vadjust{\goodbreak} a Bessel process of dimension $\delta$ which hits the
decreasing curved boundary $f(t)$ given by \eqref{eqdefde2}.

\noindent\rule{\textwidth}{0.5pt}\hypertarget{algA4}\\
\textbf{Algorithm (A4)---the square root boundary:} $\bolds{l(t)=\sqrt{\beta
_0-\beta_1 t}}$ \textbf{with} $\bolds{\beta_0>0}$, $\bolds{\beta_1>0}$.\\
Let $\kappa\in]0,1[$.\\
\textbf{Initialization:} Set $\chi(0)=0$, $\xi_0=0$, $\Xi_0=0$,
$A_0=\kappa F_\nu(\beta_0,\beta_1)$.\\
\textbf{The $\bolds{(n+1)}$th step:}
While the condition
\[
l(\Xi_{n})-\sqrt{\chi(n)}>\varepsilon\qquad\bigl(\mbox{denoted by }
\mathbb {C}(n) \bigr)
\]
holds, we define
%
\begin{equation}
\label{eqdefAn} A_n=\kappa F_\nu \biggl(\bigl(l(
\Xi_n)-\sqrt{\chi(n)}\bigr)^2, \beta_1
\biggl(1-\frac
{\sqrt{\chi(n)}}{l(\Xi_n)} \biggr) \biggr),
\end{equation}
where $F_\nu$ is defined by \eqref{eqdefdeF},
and we simulate $U_{n+1}$, a uniformly distributed random vector on
$[0,1]^{\lfloor\nu\rfloor+2}$, $G_{n+1}$, a standard Gaussian random
variable and $V_{n+1}$, a uniformly distributed random vector on
$\mathcal{S}^\delta$.
$U_{n+1}$, $G_{n+1}$ and $V_{n+1}$ have to be independent. We then
construct $(\xi_{n+1},\Xi_{n+1},\chi(n+1))$ using~\eqref{eqdefalgo}. \\
The algorithm stops when $\mathbb{C}(n)$ is not longer satisfied: we
set $\xi_{n+1}=0$ and so $\Xi_{n+1}=\Xi_{n}$ and $\chi(n+1)=\chi(n)$.\\
\rule{\textwidth}{0.5pt}

\begin{proposition}
\label{propcurved}
The statements of Theorems \ref{thmalgoline} and \ref
{thmalgo-convline} are true for Algorithm \textup{(\hyperlink{algA4}{A4})} associated with
the square root boundary.
\end{proposition}
\begin{pf}
All the arguments developed for decreasing boundaries with
lower-bounded derivatives are easily adapted to the square root
boundary. We leave the details to the reader and focus our attention to
the following fact: the stochastic process $(\Xi_n,\chi(n),  n\ge0)$
stays in the domain $D^0$ defined by
\[
D^0=\bigl\{ (t,x)\in\mathbb{R}_+^2\dvtx l(t)-\sqrt{x}>0
\bigr\}.
\]
In the $\mathit{WoMS}$ setting, for $t\ge s$, this means that for $(\Xi_n,\chi
(n))=(s,x)\in D^0 $ the following step leads to $\sqrt{\chi(n+1)}<l(\Xi
_{n+1})$. By \eqref{eqdefalgo}, it suffices to prove that
%
\begin{eqnarray}
\label{eqa-dem} \sqrt{x}+\psi_{A}(t)<l(s+t)
\nonumber
\\[-8pt]
\\[-8pt]
\eqntext{\mbox{for all } t\in\bigl
\{u\ge0\dvtx \min \bigl(l(s+u),\psi_A(u)\bigr)\ge0\bigr\},}
\end{eqnarray}
with $ A=\kappa F_\nu ( (l(s)-\sqrt{x})^2,\beta_1 (1-\frac{\sqrt {x}}{l(s)} )  )$, since $\chi(n+1)\le(\sqrt{\chi(n)}+ \psi
_{A_n}(\xi_{n+1}))^2$. By Lemma \ref{lemcomparai} and due to the
coefficient $\kappa$, we have
\[
\psi_A(t)<\sqrt{\bigl(l(s)-\sqrt{x}\bigr)^2-
\beta_1 \biggl(1-\frac{\sqrt
{x}}{l(s)} \biggr) t }.
\]
Hence
\begin{eqnarray*}
&&\bigl(l(s+t)-\sqrt{x}\bigr)^2-\psi_A(t)^2\\
&&\qquad>
\bigl(\sqrt{l(s)^2-\beta_1 t}-\sqrt {x}
\bigr)^2-\bigl(l(s)-\sqrt{x}\bigr)^2+\beta_1
\biggl(1-\frac{\sqrt{x}}{l(s)} \biggr)t
\\
&&\qquad> 2\sqrt{x} \bigl( l(s)-\sqrt{l^2(s)-\beta_1 t} \bigr)-
\frac{\beta_1\sqrt {x}}{l(s)} t\ge0.
\end{eqnarray*}
This leads directly to \eqref{eqa-dem}.
\end{pf}
%
\begin{rem}
The whole study points out a new efficient algorithm in order to
simulate Bessel hitting times for given levels or curved boundaries. We
can use this algorithm in two generalized situations:
\begin{longlist}[(1)]
\item[(1)] We have assumed that the Bessel process starts from the
origin. Of course the procedure presented here can also be applied to
Bessel processes starting from $x>0$. It suffices to change the
initialization step!
\item[(2)] We focused our attention to the Bessel process, but we
linked also the initial problem to the exit time of a $\delta
$-dimensional Brownian motion from a ball of radius~$l$. Algorithm (\hyperlink{algA1}{A1})
extended to higher dimensions can also be used in order to evaluate
exit times of general compact domains whose boundary is regular.
\end{longlist}

\end{rem}
%


\section{Numerical results}
In this part we will illustrate the previous results on some numerical
examples. Let us figure first an outcome of our algorithm, the exit
position from a sphere with radius depending on time. The figure
below is giving this result for an radius $l=1$ and a precision
$\varepsilon= 10^{-3}$.

Let us compare our algorithm with existing results. Consider the
classical Euler scheme for a Brownian motion, and evaluate the first
hitting time and hitting position from a disk with given radius.\vspace*{9pt}

\includegraphics{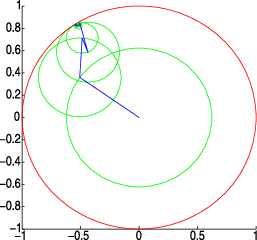}

%

\begin{figure}

\includegraphics{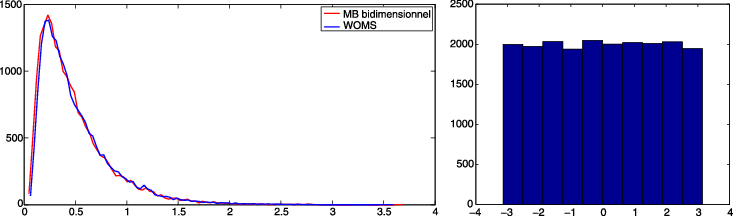}

\caption{Distribution of the hitting time [Euler scheme and $\mathit{WoMS}$
Algorithm~\textup{(\protect\hyperlink{algA1}{A1})}]---Histogram of the angle for the exit position.}
\label{fighistopossortie}\vspace*{3pt}
\end{figure}

First of all we can verify that the distribution of the hitting time
for the $\mathit{WoMS}$ algorithm matches the distribution of the hitting time
of a given level for the $2$-dimensional Bessel process. Figure \ref
{fighistopossortie} gives this result for a starting disk with
radius~$1$, a precision $\varepsilon=10^{-3}$ and a number of
simulations $N=20\mbox{,}000$. In the Euler scheme the time step is $\Delta t =
10^{-4}$.
%

We can also test the fact that the exit position is uniformly
distributed on the circle.
In order to do this we can evaluate the angle of the exit position in
our $\mathit{WoMS}$ procedure and show that it
is a uniformly distributed random variable with values in $[-\pi, \pi
]$. Figure
\ref{fighistopossortie} also shows the histogram of the result for a
disk of radius $1$ an $\varepsilon=10^{-3}$
and 20,000 simulations.

\begin{figure}[b]\vspace*{3pt}
\tabcolsep=0pt
\begin{tabular}{@{}c@{\hspace*{8pt}}c@{}}

\includegraphics{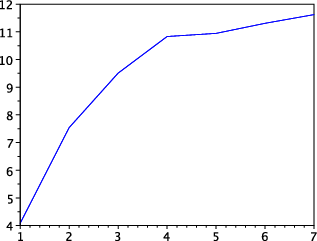}

&
\fontsize{10}{12}{\selectfont{
\begin{tabular*}{194pt}{@{\extracolsep{\fill}}lcccc@{}}
\hline
$\varepsilon$ & $10^{-1}$ & $10^{-2}$ & $10^{-3}$ & $10^{-4}$ \\
$\mathbb{E}[N^\varepsilon]$ & 4.0807 & 7.53902 & 9.50845 & 10.83133
\\[3pt]
$\varepsilon$ & $10^{-5}$ & $10^{-6}$ & $10^{-7}$ & \\
$\mathbb{E}[N^\varepsilon]$ & 10.94468 & 11.30869 & 11.62303 & \\
\hline\vspace*{100pt}
\end{tabular*}}}
\end{tabular}\vspace*{-80pt}
\caption{Averaged number of step of Algorithm \textup{(\protect\hyperlink{algA2}{A2})} versus $\varepsilon$.}
\label{graph1}
\end{figure}

Let us now present a simulation with Algorithm (\hyperlink{algA2}{A2}). We consider the
hitting time of the level $l=2$ for the Bessel process of index $\nu
=2$, and we illustrate Theorem \ref{thmalgo} by Figure \ref{graph1}.
The curve represents the averaged number of steps versus the precision
$\varepsilon=10^{-k}$, $k=1,\ldots, 7$. We can observe that the number
of steps is better than suspected since the curve is sub-linear. We
obtain the following values (for $\gamma=0.9$ and $100\mbox{,}  000$
simulations in order to evaluate the mean).

\begin{figure}
\tabcolsep=0pt
\begin{tabular}{@{}c@{\hspace*{4pt}}c@{}}
\fontsize{10}{12}{\selectfont{
\begin{tabular*}{194pt}{@{\extracolsep{\fill}}lccccc@{}}
\hline
$\nu$ & 0 & 0.5 & 1 & 1.5 & 2 \\
$\mathbb{E}[N^\varepsilon]$ & 6.819 & 7.405 & 8.270 & 8.887 & 9.594
\\[3pt]
$\nu$ & 2.5 & 3 & 3.5 & 4 & \\
$\mathbb{E}[N^\varepsilon]$ & 10.256 & 10.542 & 10.995 & 11.096 & \\
\hline\vspace*{100pt}
\end{tabular*}}}
&

\includegraphics{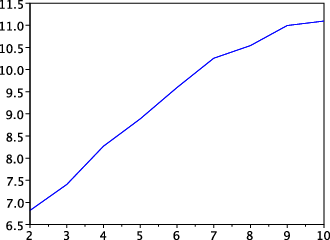}

\end{tabular}\vspace*{-80pt}
\caption{Averaged number of step of Algorithm \textup{(\protect\hyperlink{algA2}{A2})} versus $\delta=2\nu+2$.}\label{fig4}
\end{figure}

Finally we present the dependence of the averaged number of steps of
Algorithm~(\hyperlink{algA2}{A2}) with respect to the dimension of the Bessel process.
See Figure~\ref{fig4}. For
that purpose, we simulate hitting time of the level $l=2$ with
$\varepsilon=10^{-3}$, $\gamma=0.9$, $50\mbox{,}  000$ simulations for each
estimation of the averaged value, and the dimension of the Bessel
process takes value in the set $\{2,3,\ldots,18\}$.

\section{Application to the Cox--Ingersoll--Ross process}
We now aim to estimate the hitting time of a level
$l>0$ for $(X^\delta_t,  t\geq0)$, a Cox--Ingersoll--Ross process.
The CIR process is the solution of the following stochastic
differential equation:
%
\begin{equation}
\label{CIRdef} \cases{ %
\dint X^\delta_t=
\bigl(a+bX^\delta_t\bigr)\,\dint t+c\sqrt{\bigl|X^\delta_t\bigr|}
\,\dint B_t,
\vspace*{2pt}\cr
X^\delta_0=x_0,}
\end{equation}
where $x_0\geq 0$, $a\geq 0$, $b\in\mathbb{R}$, $c>0$ and $(B_t,
t\geq 0)$ is a standard Brownian motion. We denote here $\delta=
4a/c^2$.

We will first recall a connection between this stochastic process and
$(Y^\delta(t),  t\geq 0)$, the square of the Bessel process
BESQ($\delta$), the solution of the equation
%
\begin{equation}
\label{defbessel} Y^\delta(t)=y_0+\delta t+2\int
_0^t \sqrt{\bigl|Y^\delta(s)\bigr|}\,\dint
B_s,\qquad t\geq 0.
\end{equation}

\begin{lemma}
\label{lemidentite}
The CIR process has the same distribution as $(\overline{X}_t,  t\geq
0)$ which is defined by
%
\begin{equation}
\label{eqcir-relation} \cases{ %
\displaystyle\overline{X}_t=e^{bt}Y^\delta
\biggl( \frac{c^2}{4b}\bigl(1-e^{-bt}\bigr) \biggr),
\vspace*{2pt}\cr
\displaystyle\overline{X}_0=Y^\delta(0), }
\end{equation}
where $Y$ is the square of a Bessel process in dimension $\delta
=4a/c^2$; see \cite{Revuz-Yor-99}.
\end{lemma}
\begin{pf} Let us only sketch some ideas of the proof.
Let $Y^\delta(t)$ be the square of the $\delta$-dimensional Bessel
process. By applying It{\^o}'s formula, we get the stochastic
differential equation satisfied by the process $\overline{X}_t$,
%
\begin{eqnarray}
\label{eqeds1} &&\dint\overline{X}_t=b\overline{X}_t \,\dint
t+e^{bt} \,\dint \biggl(Y \biggl( \frac{c^2}{4b}
\bigl(1-e^{-bt}\bigr) \biggr) \biggr)
\nonumber
\\
&&\qquad=b\overline{X}_t \,\dint t+b\delta\frac{c^2}{4b} \,\dint
t+2e^{bt}\sqrt {\bigl|e^{-bt}\overline{X}_t\bigr|} \,\dint
B_{ ({c^2}/{(4b)})(1-e^{-bt})}
\\
&&\qquad=(a+b\overline{X}_t) \,\dint t+2e^{{bt}/{2}}\sqrt{|\overline
{X}_t|} \,\dint B_{ ({c^2}/{(4b)})(1-e^{-bt})},\nonumber
\end{eqnarray}
where $\delta=4a/ c^2.$
Let us remark that
\[
\frac{c^2}{4b}\bigl(1-e^{-bt}\bigr)=\int_0^t
\rho^2(s)\,\dint s \qquad\mbox{with } \rho(t)=\frac{c}{2}
e^{-{bt}/{2}}.
\]
We can deduce that there exists a Brownian motion
$(\beta_t,  t\geq 0)$ such that
\[
B_{ ({c^2}/{(4b)})(1-e^{-bt})}=\int_0^t \rho(s)\,\dint
\beta_s
\]
for all $t\geq 0$. With this notation, equation \eqref{eqeds1} gives
\begin{eqnarray*}
\dint\overline{X}_t&=& (a+b\overline{X}_t) \,\dint
t+2e^{{bt}/{2}}\sqrt{|\overline{X}_t|}\rho(t) \,\dint
\beta_t
\\
&=&(a+b\overline{X}_t) \,\dint t+c\sqrt{|\overline{X}_t|}
\,\dint\beta_t,
\end{eqnarray*}
and $\overline{X}_0=Y(0)$. This proves that the process $(\overline
{X}_t,  t\geq0)$ has the same distribution as the CIR process given by
(\ref{CIRdef}).
\end{pf}
Let us consider the hitting time of a given level $l$ for the CIR
process and denote it by $T_l$. This time is defined by
\[
T_l=\inf\bigl\{ s\geq 0; X^\delta_s=l\bigr\}.
\]
The previous Lemma \ref{lemidentite} gives also an equivalence (in
distribution) connecting the hitting time of the CIR process and the
hitting time of the square of a $\delta$-dimensional Bessel process.
%
\begin{proposition}\label{propenfin}
The hitting time $T_l$ of a level $l>0$ for a CIR process has the same
distribution as $-\frac{1}{b}\log (1-\frac{4b}{c^2}\tau_\psi )$ where
\[
\tau_\psi=\inf \biggl\{ t\geq 0; Y^\delta(t)=l \biggl(1-
\frac{4b}{c^2}t \biggr) \biggr\},
\]
and $Y^\delta$ is the square of a Bessel process of dimension $\delta=4a/c^2$.
\end{proposition}
\begin{pf}
By using Lemma \ref{lemidentite}, $\tau_\psi$ has the same
distribution as $\overline{T}_l$ given by
%
\begin{equation}
\label{deftau2} \overline{T}_l=\inf \biggl\{s\geq 0;
Y^\delta \biggl( \frac
{c^2}{4b}\bigl(1-e^{-bs}\bigr)
\biggr)=le^{-bs} \biggr\}.
\end{equation}
Define $t=\frac{c^2}{4b}(1-e^{-bs})$, so we have two situations:

\textit{First case}: If $b<0$, let $s=\eta(t)$ where
\[
\eta(t)=-\frac{1}{b}\log \biggl(1-\frac{4b}{c^2}t \biggr)\qquad
\mbox{for } t\geq 0.
\]
The map $\eta$ is a strictly nondecreasing function, and we thus get thus
\begin{eqnarray*}
\overline{T}_l&=&\inf \biggl\{\eta(t); t\geq 0, Y^\delta(t)=l
\biggl(1-\frac
{4b}{c^2}t \biggr) \biggr\}\\
&=&\eta \biggl( \inf \biggl\{t\geq
0;Y^\delta(t)=l \biggl(1-\frac{4b}{c^2}t \biggr) \biggr\} \biggr).
\end{eqnarray*}

\textit{Second case}: If $b\geq 0$, let also $s=\eta(t)$. In this case
the variable $t$ takes its values only on the interval $ [0,\frac
{c^2}{4b} )$. So
\[
\overline{T}_l=\inf \biggl\{\eta(t); 0\le t\le\frac{c^2}{4b},
Y^\delta (t)=l \biggl(1-\frac{4b}{c^2} t \biggr) \biggr\}.
\]
The condition $0\le t\le\frac{c^2}{4b}$ can be omitted in the
estimation of the infimum as the boundary to hit: $1-\frac{4bt}{c^2}$
is negative outside this interval, and the Bessel process is always positive.
Furthermore the function $\eta$ is also nondecreasing for $b\geq 0$,
and the result is thus obtained.
\end{pf}

\textit{Application of Algorithm} (\hyperlink{algA4}{A4}):

An immediate consequence of Proposition \ref{propenfin} is that the
hitting time $T_l$ is related to the first time the Bessel process of
dimension $\delta=4a/c^2$ reaches the curved boundary: $f(t)=\sqrt{l
(1-\frac{4b}{c^2}t  )}$. We are able to apply Algorithm (\hyperlink{algA4}{A4}) if
$4a/c^2\in\mathbb{N}^*$ and $b>0$ (the boundary is then decreasing).
Let us denote by $N^\varepsilon$ the number of steps of (\hyperlink{algA4}{A4}) and $\Xi
_{N^\varepsilon}$, the approximated hitting time of the Bessel process
associated with the particular curved boundary $f$. Combining
Propositions \ref{propcurved} and \ref{propenfin} leads to
\begin{eqnarray*}
\biggl(1-\frac{\varepsilon}{\sqrt{2\alpha\pi}} \biggr)\mathbb{P} \biggl(\Xi
_{N^\varepsilon}\le
\frac{c^2}{4b} \bigl(1-e^{-bt}\bigr)-\alpha \biggr)&\le&\mathbb
{P}(T_l\le t)\\
&\le&\mathbb{P} \biggl(\Xi_{N^\varepsilon}\le
\frac{c^2}{4b} \bigl(1-e^{-bt}\bigr) \biggr)
\end{eqnarray*}
for $\alpha$ small enough and $t>0$.
%
%
\begin{appendix}\label{app}
\section*{Appendix: Simulation of random variables}
Let us introduce simulation procedures related to particular
probability density functions.
%
\begin{propositionn}
\label{propappend}
Let $Z$ be a random variable with Gamma distribution ${\rm Gamma}(\alpha
,\beta)$, that is,
\[
\mathbb{P}(Z\in\dint z)=\frac{1}{\Gamma(\alpha)\beta^\alpha} z^{\alpha-1}e^{-{z}/{\beta}}
\ind{\{z>0 \}} \,\dint z,\qquad \alpha >0, \beta>0.
\]
Then $W=\exp(-Z)$ has the following distribution:
\[
\mathbb{P}(W\in\dint r)=\frac{1}{\Gamma(\alpha)\beta^\alpha}(-\log r)^{\alpha-1}r^{1/\beta-1}
\ind{[0,1]}(r) \,\dint r.
\]
In particular the stopping time $\tau_\psi$ defined by \eqref
{densitetau} has the same law as $ [\frac{a}{\Gamma(\nu+1)2^\nu
} ]^{{1}/{(\nu+1)}}e^{-Z}$. Here $Z$ is a Gamma distributed
random variable with parameters $\alpha=\nu+2$ and $\beta=\frac{1}{\nu+1}$.
\end{propositionn}
\begin{pf} Let $f$ be a nonnegative function. Using suitable
changes of variables, we obtain
\begin{eqnarray*}
\mathbb{E}\bigl[f(W)\bigr]&=&\frac{1}{\Gamma(\alpha)\beta^\alpha}\int_0^\infty
f\bigl(e^{-z}\bigr) z^{\alpha-1}e^{-{z}/{\beta}} \,\dint z
\\
&=&\frac{1}{\Gamma(\alpha)\beta^\alpha}\int_0^1f(r) (-\log
r)^{\alpha
-1}r^{1/\beta-1} \,\dint r.
\end{eqnarray*}
In order to end the proof it suffices to multiply $W$ by a constant and
use once again a change of variables formula.
\end{pf}
We need to simulate Gamma distributed variables. Let us just recall
some common facts.
%
\begin{propositionn}
\textup{(i)} If $\alpha\in\mathbb{N}$ (so-called Erlang distributions),
then the Gamma distributed variables $Z$ has the same law as
\[
-\beta\log(U_1\cdots U_\alpha),
\]
where $(U_i)_{1\le i\le\alpha}$ are independent uniformly distributed
random variables. Hence $W$ defined by $W=\exp(-Z)$ can be simulated by
\[
(U_1U_2\cdots U_\alpha)^\beta.\vspace*{-12pt}
\]
\begin{longlist}[(ii)]
\item[(ii)] If $\alpha-1/2\in\mathbb{N}$, then $Z$ has the same law as
\[
-\beta\log(U_1\cdots U_{\lfloor\alpha\rfloor})+\frac{\beta N^2}{2},
\]
where $(U_i)_{1\le i\le\lfloor\alpha\rfloor}$ are i.i.d. uniformly
distributed random variables, and $N$ is an independent standard
Gaussian r.v.; see, for instance, \cite{Devroye}, Chapter \textup{IX.3}.
\end{longlist}
\end{propositionn}
\end{appendix}

%



\printaddresses

\end{document}